\documentclass{amsart}
\usepackage[utf8]{inputenc}
\usepackage[margin=2cm,includeheadfoot,centering]{geometry} 
\usepackage{amsmath, amsthm, amsfonts, amssymb, xy,  mathrsfs, graphicx}
\usepackage{hyperref}
\usepackage{array}
\usepackage{cite}
\usepackage{tikz}
\usepackage[section]{placeins}
\usetikzlibrary{mindmap, backgrounds}

\newtheorem{theorem}{Theorem}[section]

\newtheorem{defi}[theorem]{Definition}
\theoremstyle{definition}

\newtheorem{remark}[theorem]{Remark}

\theoremstyle{remark}

\newtheorem{criterion}{Criterion}

\newcommand{\la}{\leftarrow}
\newcommand{\ra}{\rightarrow}
\newcolumntype{L}{>$l<$}
\newcolumntype{P}{>$p{8cm}<$}
\newcolumntype{Q}{>$p{10cm}<$}
\DeclareMathOperator{\Pf}{Pf}

\title{Minuscule Schubert varieties of exceptional type}
\author{Sara Angela Filippini, Jacinta Torres and Jerzy Weyman}
\date{\today}

\begin{document}

\begin{abstract}
We study exceptional minuscule Schubert varieties and provide the defining equations of the defining ideals of their intersection with the big open subset. We also provide the resolutions of these ideals and characterize some of them in terms of fundamental examples of ideals in the theory of Gorenstein ideals.
\end{abstract}

\maketitle
\section{Introduction}

It is a long-standing open problem to classify all Gorenstein Schubert varieties. There has been much work on this problem (\cite{brownlakshmibai,WooYong2006,Perrin2009}); however an exhaustive description does not exist. In this work, we give an open subset $U \subset \mathbb{P}(V_{\omega_{i}})$ such that the intersections $Y_{\sigma} = X_{\sigma} \cap U$, where $X_{\sigma}$ is an exceptional minuscule Schubert variety have descriptions in terms of well-known perfect ideals of low codimension. Some of the ideals we uncover have appeared in commutative algebra studies, such as the defining ideals of varieties of complexes studied by Herzog in \cite{Herzog1974} and more generally by Kustin in \cite{kustin1993}. Explicitly, we have

\begin{theorem}
\label{e6p1gorenstein}
Let $G$ be a reductive group of exceptional type and $P \subset G$ a standard parabolic subgroup stabilizing a minuscule fundamental weight $\omega_{i}$. Then there exists an open subset $U\subset  \mathbb{P}(V_{\omega_{i}})$ such that for any given minuscule Schubert variety $X_{\sigma} \subset G/P$, the intersection $Y_{\alpha} = X_{\alpha } \cap U$ is as described in Sections 3 and 4. In particular, for $\sigma$ of type $E_{6}$ or as described in Section 4, $Y_{\sigma}$  is a complete intersection in one of the following:
\begin{enumerate}
\item[a.] $Y_{\sigma}$ is a complete intersection -  the minimal free resolution of its coordinate ring is a Koszul complex.
\item[b.] in the codimension three variety of submaximal Pfaffians of a skew symmetric matrix.
\item[c.]  in the variety of pure spinors. 
\item[d.]  in a variety of complexes. 
\item[e.] in a Huneke-Ulrich ideal of deviation 2. 
\item[f.] in the variety defined by the vanishing of the $2 \times 2$ minors of a $2\times 3$ generic matrix. 
\item[g.] in the variety defined by the vanishing of $4 \times 4$ Pfaffians of a $6\times 6$ skew-symmetric matrix. \end{enumerate}
\end{theorem}

The paper is organized in three sections: Section 2 is dedicated to preliminaries on Schubert varieties, linkage and deformation theory. We will collect all results we need to carry out our study in Sections 3 and 4. Section 3 is dedicated to $E_{6}$ and Section 4 to $E_{7}$. This is due to the fact that there are no minuscule fundamental weights of type $E_{8}$. In both sections we first determine the defining equations for all Schubert varieties $X_{\sigma} \subset G/P$ using the work of Vavilov, Luzgarev and Pevzner \cite{vavilove6}, respectively \cite{vavilove7} and compute their Hilbert polynomials using Macaulay2. Since Schubert varieties are Cohen Macaulay we may apply a result of Stanley \cite{palindromic} to check which of these are Gorenstein, simply by checking which of the corresponding Hilbert polynomials are palindromic. We equal our list with the one obtained by Brown-Lakshmibai in \cite{brownlakshmibai}.  Another motivation for our work is the forthcoming work on the structure of Gorenstein ideals of low codimension, as several key examples turn out to be open subsets of Schubert varieties in homogeneous spaces for exceptional groups (see for example $J_{51}$ in Section 4 and $I_{23}$ in Section 3).  For analogous examples of perfect ideals of codimension 3, see \cite{sam2020schubert}.   \\

 To prove Theorem \ref{e6p1gorenstein}, we manipulate the given equations accordingly by restricting the fundamental representation corresponding to $P$ (which we assume to be the 27-dimensional representation $V(\omega_{1})$ without loss of generality) to the Levi subgroup of $G$ of type $D_{5}$. We study the remaining Schubert varieties $X_{\sigma}$ individually and calculate minimal free resolutions of their defining ideals as well as Hilbert functions using Macaulay2. Moreover, by a geometric method we calculate the expected terms of a resolution $R_{\sigma}$ which  ``is" essentially the resolution of the coordinate ring of $X_{\sigma}$. We obtain generators of an ideal whose associated resolution is precisely $R_{\sigma}$, and construct our resolution in terms of $R_{\sigma}$. We organize our results in various tables, where we also include other data where available, including: symmetry group, reduced expression of $\sigma$ in $W/W_{P}$ in terms of simple reflections, and quick description of $Y_{\sigma}$. For the Gorenstein ideals of codimension four, we recover the explicit descriptions by Celikbas-Laxmi-Weyman and the spinor structures on their resolutions.

\section{Basics on Schubert varieties and commutative algebra} 
\subsection{Linkage and residual intersections}
\begin{defi}[\cite{HunekeUlrich}]
Let $R$ be a commutative Noetherian ring and $I,J \subset R$ two ideals. Then 

\[I:J = \left\{f \in R : f*J \subset I \right\} \]

\noindent
is an ideal, called the colon ideal of $I$ and $J$. If $I$ is generated by elements $a_{1},..., a_{s} \in J$ then if $\operatorname{ht}(J) \leq s \leq \operatorname{ht}(A)$ we say that $A = I:J$ is a \textbf{residual intersection} of $J$. If the sequence $(a_{1},..., a_{s})$ is regular, we say that $I$ and $J$ are \textbf{linked}.
\end{defi}

\subsection{Complete intersections}

\begin{defi}
Let $R$ be a Noetherian commutative ring. An ideal $I$ is a \textit{complete intersection} ideal of codimension $n$ if there exists a regular sequence of length $n (a_{1},...,a_{n})$ such that $I = (a_{1},...,a_{n})$. 
\end{defi}

\begin{theorem}\cite[Corollary 1.2.3 ]{WeymanBook}
\label{ci}
Let $R$ be a Noetherian commutative ring and $I = (a_{1},...,a_{n})$ a complete intersection ideal. Then the Koszul complex on the elements $a_{1},...,a_{n}$ is a free resolution of $R/I$. 
\end{theorem}

\begin{defi}
Let $(d_{1},...,d_{n})$ be a sequence of positive integers, and $(a_{1},...,a_{n})$ a regular sequence in a polynomial ring $k[x_{1},...,x_{s}]$ such that $\operatorname{deg}(a_{i}) = d_{i}$. Let $J = (a_{1},...,a_{n})$. If there exist varieties $Y' \subset Y$ and an exact sequence
\[0 \rightarrow J \rightarrow k[Y] \twoheadrightarrow k[Y']\]

then we say that $Y'$ is a complete intersection of type $(d_{1},...,d_{n})$ in $Y$. 
\end{defi}

For some of the analysis we are carrying out about the Schubert varieties, we will need the following criterion from \cite{HunekeUlrich1987} (see p. 279 and Corollary 5.13):

\begin{criterion} [Licci criterion]
\label{licci}
Let $R = k[x_{1},...,x_{n}]$ and $I \subset R$ be a homogeneous ideal with minimal graded free resolution,

\[ 0 \rightarrow \overset{b_g}{\underset{j = 1}{\bigoplus}} R(-n_{gj}) \rightarrow \cdots \rightarrow  \overset{b_1}{\underset{j = 1}{\bigoplus}} R(-n_{1j}) \rightarrow R \rightarrow R/I \rightarrow 0\]

\noindent 
where $g$ is the height (or codimension) of $I$. If 

\begin{align}
\label{crit}
\operatorname{max}\left\{n_{gj}\right\} \leq (g-1)\operatorname{min}\left\{n_{1j}\right\}
\end{align}

\noindent
 then $(R/I)_{(x_{1},...,x_{n})}$ is not in the linkage class of a complete intersection. 
\end{criterion}

\subsection{Schubert varieties}
Let  $G$ be a reductive algebraic group and let $P_{\omega_{i}} \subset G$ be a parabolic subgroup which stabilizes a fundamental weight $\omega_{i}$, using Bourbaki notation to label our fundamental weights. It is known  that $ G/P_{\omega_{i}}$ is naturally embedded into $\mathbb{P}(V_{\omega_{i}})$, where $V(\omega_{i})$ is a fundamental representation of highest weight $\omega_{i}$. Let $B\subset P_{\omega_{i}}$  be a Borel subgroup. The Bruhat decomposition

\[G/P_{\omega_{i}} = \underset{w \in W/W_{P_{\omega_{i}}}}{\bigsqcup} B \cdot w P_{\omega_{i}} \]

\noindent 
implies that the Schubert varieties defined as the orbit closures 
\[X_{w} : = \overline{B \cdot w} \subset G/P_{\omega_{n}}\]

 have the following cell decomposition: 

\begin{align}
\label{celldec}
X_{w} = \underset{\underset{w \in W/W_{P_{\omega_{n}}}}{v \leq w}}{\bigsqcup} B v. 
\end{align}

The quotient $W/W_{P_{\omega_{n}}}$ is a Coxeter group with generators $\sigma_{1},...,\sigma_{n}$, each of which is the projection into the quotient of the corresponding simple reflection in $W$. The Bruhat graph associated to the fundamental weight $\omega_{i}$ is the directed colored graph whose vertices are the elements of $W/W_{P_{\omega_{i}}}$ and such that there exists an arrow $v \overset{t}{\rightarrow} w$ if and only if, for minimal length representatives $v',w' \in W, w' = v's$ for the simple reflection $s_{t}$. The relation $v \leq w$ in (\ref{celldec}) holds if and only if there exists a minimal path $v \overset{t_{1}}{\rightarrow} \cdots \overset{t_{k}}{\rightarrow} w$ in the Bruhat graph. If the fundamental weight $\omega_{i}$ is \textbf{minuscule}, this graph coincides with the so-called \textbf{crystal graph} associated to the fundamental representation $V_{\omega_{i}}$.
%
%


\section{Type $E_{6}$}
\subsection{Minuscule homogeneous spaces}

Let $G$ be a simple algebraic group of type $E_6$ and $P \subset G$ a parabolic subgroup which stabilizes a minuscule fundamental weight $\omega_{i}$. If we use Bourbaki notation to label our fundamental weights, we have that the only minuscule fundamental weights in this case are $\omega_{1}$ and $\omega_{6}$. 

It is known  that $G/P_{i}$ is naturally embedded into $\mathbb{P}(V_{\omega_{i}})$, where $V(\omega_{i})$ is the 27-dimensional fundamental representation  (for $i = 1$ or $6$). Since $V(\omega_{1})$ and $V(\omega_{6})$ are dual to each other, we can assume $i =1$.  In fact, we know that $G/P = \overline{G\cdot v_{\omega_{1}}}$, where $v_{\omega_{1}}$ is the highst weight vector in $V(\omega_{1})$. We identify the coordinate ring of $V(\omega_{1})$ with the polynomial ring on $27$ variables $x_{1},...,x_{27}$. Let us denote this ring by $R = k[x_{1},...,x_{27}]$. 

\subsection{Defining equations}
In \cite{vavilove6}, Vavilov-Luzgarev-Pevzner have obtained as defining equations for $G/P$ the $27$ partial derivatives of the following cubic invariant:  

\begin{align}
\label{cubic}
Q &= x_{1}x_{18}x_{27} - x_{1}x_{19}x_{26} + x_{1}x_{20}x_{25} -  x_{1}x_{21}x_{24}+  x_{1}x_{22}x_{23} \nonumber \\ 
&-x_{2}x_{11}x_{27} + x_{2}x_{13}x_{26 }- x_{2}x_{15}x_{25}  + x_{2}x_{16}x_{24} - x_{2}x_{17}x_{23} \nonumber \\
&+x_{3}x_{9}x_{27} - x_{3}x_{12}x_{26} +x_{3}x_{14}x_{25} - x_{3}x_{16}x_{22} + x_{3}x_{17}x_{21} \nonumber \\
&-x_{4}x_{7}x_{27} +x_{4}x_{10}x_{26} -x_{4}x_{14}x_{24} +x_{4}x_{15}x_{22} -x_{4}x_{17}x_{20} \nonumber \\
&+x_{5}x_{6}x_{27} -x_{5}x_{8}x_{26} + x_{5}x_{14}x_{23} -x_{5}x_{15}x_{21} +x_{5}x_{16}x_{20} \\
&-x_{6}x_{10}x_{25} + x_{6}x_{12}x_{24} - x_{6}x_{13}x_{22} + x_{6}x_{17}x_{19} + x_{7}x_{8}x_{25} \nonumber \\
&-x_{7}x_{12}x_{23} + x_{7}x_{13}x_{21} - x_{7}x_{16}x_{19} - x_{8}x_{9}x_{24} + x_{8}x_{11}x_{22}  \nonumber \\
&-x_{8}x_{17}x_{18} + x_{9}x_{10}x_{23}  - x_{9}x_{13}x_{20} + x_{9}x_{15}x_{19} - x_{10}x_{11}x_{21} \nonumber \\
&+x_{10}x_{16}x_{18} + x_{11}x_{12}x_{20} - x_{11}x_{14}x_{19} - x_{12}x_{15}x_{18} + x_{13}x_{14}x_{18} \nonumber
\end{align}
For convenience the equations for the cubic and its derivatives are listed in the appendix \ref{sec:appendix}
That is, 

\begin{align}
\label{orbitdescription}
G/P = \left\{ [v] \in \mathbb{P}(V(\omega_{1})) |  \frac{\partial Q }{\partial x_{i}}(v) = 0 \hbox{    } \forall i \in [27] \right\}
\end{align}

We are going to consider the restriction of $V(\omega_{1})$ to the Levi subgroup of type $D_{5}$ which corresponds to deleting the first node (in Bourbaki order) from the $E_{6}$ Dynkin diagram. The result of this branching is as follows (this can be easily seen on the Bruhat graph, which in this case corresponds to the crystal graph of $V(\omega_{1})$): 

\begin{align}
\operatorname{res}V(E_{6}, \omega_{1})_{D_{5}}^{E_{6}} = V(D_{5},0) \oplus V(D_{5}, \omega_{4}) \oplus V(D_{5}, \omega_{1}) 
\end{align}

Minuscule Schubert varieties are labelled by the minimal length representatives in $W/W_{P}$.  In the graph depicted below, we label the nodes by using two labelling systems: the first, with numbers from $1$ to $27$ is the ``D5 labelling'' given by Vavilov, Luzgarev and Pevzner in \cite{vavilove6}, and the second is with the variables $y_{ij}, y_{ijkl}$ for $1\leq i < j < k < l \leq 5$, which correspond to the even spinors, that is, the \textit{positive} half-spinor and fourth fundamental representation  of $\mathfrak{so}(10,\mathbb{C})$, the variables $z_{i}, \bar z_{i}$, which correspond to the basis of the natural and and first fundamental representation of $\mathfrak{so}(10,\mathbb{C})$, and the label $x$, which corresponds to the trivial representation $\mathbb{C}$.

\begin{figure}
\label{graphe6p1}
\includegraphics[scale = 0.6]{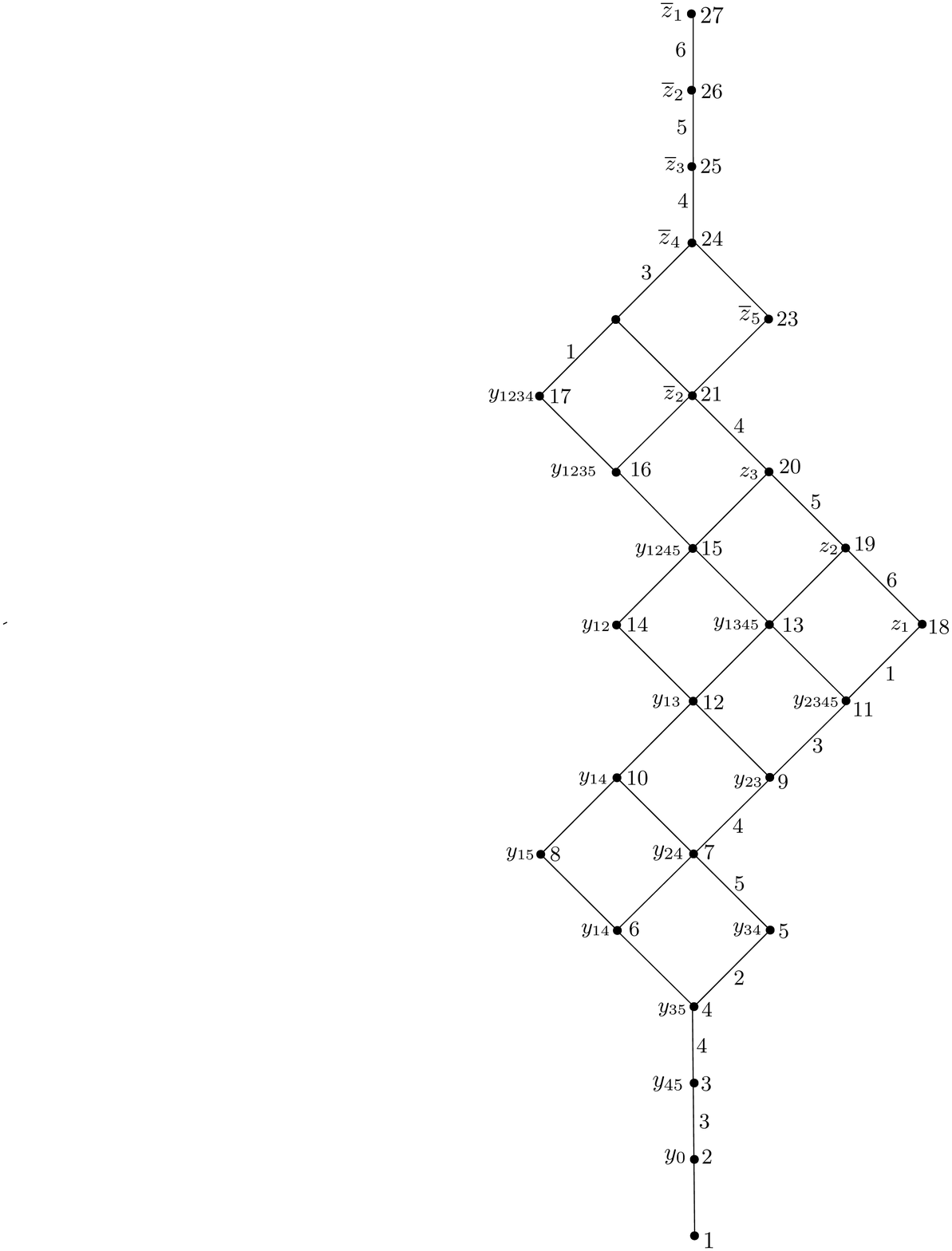}
\caption{Bruhat interval $W/W_{P_{1}}$ for $E_{6}$}
\end{figure}

Taking the 27 partial derivatives of $Q$ and substituting $x_{i}$ by the corresponding $l$ in our labelling from Figure \ref{graphe6p1}  $\left(\begin{smallmatrix} i \\ l \end{smallmatrix} \right)$, we obtain the following equations:

We will in fact work in the open subset $U_{x}\subset \mathbb{P}(V(\omega_{1}))$ defined by $x\neq 0$. The intersection $Y := U_{x} \cap G/P_{1}$ is then defined as the zero set of  the 27 equations above when substituting $x = 1$. Moreover, it is easy to see (using, for instance, Macaulay2), that in this case the vanishing of the last ten polynomials (after making the substitution $x = 1$), $f_{18},...,f_{27}$ is a minimal set of generators of the ideal defining the open cell $Y$. For each node $l$ in Figure \ref{graphe6p1}, we will denote the corresponding Schubert variety by $X_{l}$ and by $Y_{l}$ its intersection with the open cell $Y$.

\subsection{Schubert varieties in $G(E_6)/P_1$}
\subsubsection{Tables}

\

${\small\renewcommand{\arraystretch}{1.2}
\begin{array}{llllL}
\mbox{ideal} & \mbox{codim} & \mbox{resolution} & \mbox{numerator} & description \\
\hline
I_{27} & 0 & expected & 1 & open subset in $G(E_6)/P_1$ \\
\hline I_{26} & 1 & expected & 1+T & c.i. of type $(1)$, Gorenstein\\
\hline I_{25} & 2 & expected & 1+2 T+T^{2} &c. i. of type $(2,2)$, Gorenstein \\
\hline I_{24} & 3 & expected & 1+3 T+3 T^{2}+T^{3} & c.i. of type (2,2,2), Gorenstein  \\
\hline I_{23} & 4 &  & 1+3 T+T^{2} & c.i. of type (1) in ideal of submax. Pf. of $5\times 5$ skew-symm, Gorenstein \\
\hline I_{22} & 4 & & 1+4 T+5 T^{2}+T^{3} & almost c.i. \\
\hline I_{21} & 5 &  & 1+4 T+4 T^{2}+T^{3} & c.i. of type (1,2) in the ideal of submax. Pf. of $5\times 5$ skew-symm, Gorenstein\\
\hline I_{20} & 16 & & 1+4 T+3 T^{2} & c.i. of type (1,1) in a a variety of complexes \\
\hline I_{19} & 7 & & 1+3 T &  c.i. of type $(1^4)$ in EN of $2\times 4$\\
\hline I_{18} & 8 &  & 1 & non-singular  \\
\hline I_{17} & 5 & expected & 1+5 T+5 T^{2}+T^{3} &  [KW], [G], Gorenstein\\
\hline I_{16} & 6 & & 1+5 T+5 T^{2}+T^{3} & c.i. of type $(1)$ in $I_{17}$, Gorenstein \\
\hline I_{15} & 7 & & 1+5 T+5 T^{2}+T^{3} & c.i. of type $(1,1)$ in $I_{17}$, Gorenstein \\
\hline I_{14} & 8 & & 1+3 T+ T^{2} &c.i. of type $(1^5)$ in  Pf. of $5\times 5$ skew, Gorenstein  \\
\hline I_{13} & 9 & & 1+4 T+2 T^{2} & c.i. of type $(1^4)$ in a variety of complexes \\
\hline I_{12} & 9 & & 1+3 T+ T^{2} &c.i.of type $(1^6)$ in Pf. of $5\times 5$ skew , Gorenstein\\
\hline I_{11} & 9 & & 1 + T & c.i. of type $(1^8,2)$, Gorenstein \\
\hline I_{10} & 10 & &  1 + 2 T & c.i. of type $(1^8)$ in Eagon-Northcott for $2\times 3$ matrices\\
\hline I_{9} & 10 & & 1 + T & c.i. of type $(1^9,2)$, Gorenstein  \\
\hline I_{8} & 11 & & 1 & non-singular \\
\hline I_{7} & 11 &  &  1 + T & c.i. of type $(1^10,2)$, Gorenstein \\
\hline I_{6} & 12 & & 1 & non-singular \\
\hline I_{5} & 12 & & 1 & non-singular \\
\hline I_{4} & 13 & & 1 & non-singular \\
\hline I_{3} & 14 & & 1 & non-singular \\
\hline I_{2} & 15 & & 1 & non-singular \\
\hline I_{1} & 16 & & 1 & non-singular \\
\end{array}}
$


\subsection{($I_{26}$, $E_6$)}\label{i26e6} 
$I_{26} = (-{y}_{1234} {y}_{15}+{y}_{1235} {y}_{14}-{y}_{1245}
      {y}_{13}+{y}_{1345} {y}_{12})$ is resolved by Koszul complex on one quadric $R \leftarrow R(-2) \leftarrow 0$, $P(I_{26}) = 1-T^{2}$. 
            
\subsection{($I_{25}$, $E_6$)}\label{i25e6} 
$I_{25} = (-{y}_{1234} {y}_{15}+{y}_{1235} {y}_{14}-{y}_{1245}
      {y}_{13}+{y}_{1345} {y}_{12},-{y}_{1234} {y}_{25}+{y}_{1235}
      {y}_{24}-{y}_{1245} {y}_{23}+{y}_{2345} {y}_{12})$ is resolved by Koszul complex on two quadrics $R \leftarrow R^{2}(-2) \leftarrow R(-4) \leftarrow 0$, $P(I_{25}) = 1+2 T+T^{2}$.      
     
\subsection{($I_{24}$, $E_6$)}\label{i24e6} 
\begin{align*}
I_{24} = & (-{y}_{1234} {y}_{15}+{y}_{1235} {y}_{14}-{y}_{1245} 
      {y}_{13}+{y}_{1345} {y}_{12},-{y}_{1234} {y}_{25}+{y}_{1235} {y}_{24}-{y}_{1245} {y}_{23}+{y}_{2345} {y}_{12}, \\
      &-{y}_{1234}
      {y}_{35}+{y}_{1235} {y}_{34}-{y}_{1345} {y}_{23}+{y}_{2345} {y}_{13})
      \end{align*}   
is resolved by Koszul complex on three quadrics $R \leftarrow R^{3}(-2) \leftarrow R^{3}(-4) \leftarrow R(-6) \leftarrow 0$, $P(I_{24}) = -T^6 + 3T^4 - 3T^2 + 1$.


\subsection{($I_{23}$, $E_6$)}\label{i23e6} 
\begin{align*}  I_{23} = & ({y}_{1235} {y}_{14}-{y}_{1245}
      {y}_{13}+{y}_{1345} {y}_{12}, {y}_{1235}
      {y}_{24}-{y}_{1245} {y}_{23}+{y}_{2345} {y}_{12}, {y}_{1235} {y}_{34}-{y}_{1345} {y}_{23}+{y}_{2345} {y}_{13},\\
  & {y}_{1245} {y}_{34}-{y}_{1345}{y}_{24}+{y}_{2345} {y}_{14}, {y}_{14} {y}_{23}-{y}_{24} {y}_{13}+{y}_{34} {y}_{12},{y}_{1234})
      \end{align*} 
is a complete intersection of type $(1)$ in the codimension three Gorenstein ideal of submaximal Pfaffians of the $5\times 5$ skew-symmetric matrix      
$$M_{23} = \begin{pmatrix}0&
      {y}_{1235}&
      {y}_{1245}&
      {y}_{1345}&
      {y}_{2345}\\
      {-{y}_{1235}}&
      0&
      {y}_{12}&
      {y}_{13}&
      {y}_{23}\\
      {-{y}_{1245}}&
      {-{y}_{12}}&
      0&
      {y}_{14}&
      {y}_{24}\\
      {-{y}_{1345}}&
      {-{y}_{13}}&
      {-{y}_{14}}&
      0&
      {y}_{34}\\
      {-{y}_{2345}}&
      {-{y}_{23}}&
      {-{y}_{24}}&
      {-{y}_{34}}&
      0\\
      \end{pmatrix}.$$
 
\noindent
The resolution of $I_{23}$ can be obtained by tensoring the resolution of the ideal generated by the Pfaffians, $R \leftarrow R^{5}(-2) \leftarrow R^{5}(-3) \leftarrow R(-5) \leftarrow 0$, with Koszul complex in $1$ variable $R \leftarrow R(-1)$, $P(I_{23}) = (1+3T+T^2)(1-T)^3 = -T^{5} + 5T^3 - 5T^2 + 1.$   
%

\begin{remark}
This is an open cell in the only Gorenstein Schubert variety of codimension 4, and is the expected generic form of the resolution of the Gorenstein ideal of codimension 4 and 6 generators studied in  \cite{HerzogMiller1985}. 
\end{remark}

\

\subsection{($I_{22}$, $E_6$)}\label{i22e6} 
\begin{align*}   I_{22} = \, & (-{y}_{1234} {y}_{15}+{y}_{1235} {y}_{25}-{y}_{1245}
       {y}_{35}+{y}_{1235} {y}_{45},-{y}_{2345} {y}_{14}+{y}_{1345}
       {y}_{24}-{y}_{1245} {y}_{34}+{y}_{1234} {y}_{45}, -{y}_{2345} {y}_{13}+{y}_{1345} {y}_{23}
       \\
       &
       -{y}_{1235} {y}_{34}+
       {y}_{1234}
       {y}_{35},
       -{y}_{2345} {y}_{12}+{y}_{1245} {y}_{23}-{y}_{1235}
       {y}_{24}+{y}_{1234} {y}_{25},-{y}_{1345} {y}_{12}+{y}_{1245}
       {y}_{13}-{y}_{1235} {y}_{14}+{y}_{1234} {y}_{15})
  \end{align*}
  is an almost complete intersection with $5$ generators with resolution \[ R \leftarrow R^{5}(-2) \leftarrow R^{11}(-4)\oplus R(-3) \leftarrow R^{10}(-5) \leftarrow R(-7)\oplus R(-6) \leftarrow 0,  \]
  $P(I_{22}) = T^7 + T^6- 10t^5+ 11t^4 + t^3 - 5t^2 + 1$.

\begin{remark}
$I_{22}$ and $I_{23}$ are linked by the $4$ quadrics: 

\begin{align*}
 y_{15}y_{1234} - y_{25}y_{1235} + y_{13}y_{1245} - y_{1345}y_{12} \\
 y_{25}y_{1234} - y_{24}y_{1235} + y_{23}y_{1245} - y_{2345}y_{12} \\
 y_{13}y_{1234} - y_{34}y_{1235} + y_{23}y_{1345} - y_{2345}y_{13} \\
 y_{45}y_{1234} - y_{34}y_{1245} + y_{24}y_{1345} - y_{25}y_{2345}.
 \end{align*}
 
 \noindent
 This is readily checked in Macaulay2 and is an example of the work of Ulrich \cite{Ulrich1988}, where two codimension $g$, linked subvarieties intersect in a Gorenstein subvariety of codimension $g+1$. 
\end{remark}

\subsection{($I_{21}$, $E_6$)}\label{i21e6} 
\begin{align*} 
I_{21} = \, & (y_{14}y_{1235} - y_{13}y_{1245} + y_{1345}y_{12},
      y_{24}y_{1235} - y_{23}y_{1245} + y_{2345}y_{12},  y_{34}y_{1235} - y_{23}y_{1345} + y_{2345}y_{13},  \\
   &  
      y_{34}y_{1245} - y_{24}y_{1345} + y_{14}y_{2345},   y_{34}y_{12} - y_{24}y_{13} + y_{23}y_{14}, -y_{45}y_{1235} + y_{35}y_{1245} - y_{25}y_{1345} + y_{15}y_{2345} ,\\
  &   y_{1234})
\end{align*}
is a a complete intersection of type $(1,2)$ in  the codimension three ideal of submaximal Pfaffians of $M_{23}$.
Its resolution 
%
%
is given by tensoring together $R \leftarrow R^{5}(-2) \leftarrow R^{5}(-3) \leftarrow R(-5) \leftarrow 0$, the Koszul complex of the degree $1$ variable $y_{1234}$, and the one for the degree $2$ variable $-{y}_{1235} {y}_{45}+{y}_{1245}
       {y}_{35}-{y}_{1345} {y}_{25}+{y}_{2345} {y}_{15}$, $P(I_{21}) = T^7 - 6T^{5} +5T^4 + 5T^3 - 6T^2 + 1$.
%
\subsection{($I_{20}$, $E_6$)}\label{i20e6} 
$I_{20}$ generated by:
\begin{align*}
- y_{15}y_{1234} + y_{14}y_{1235} - y_{13}y_{1245} + y_{12}y_{1345}\\
 - y_{25}y_{1234} + y_{24}y_{1235} - y_{23}y_{1245} + y_{12}y_{2345}\\
 - y_{35}y_{1234} + y_{34}y_{1235} - y_{23}y_{1345} + y_{13}y_{2345}\\
 - y_{45}y_{1234} + y_{34}y_{1245} - y_{24}y_{1345} + y_{14}y_{2345}\\
 - y_{45}y_{1235} + y_{35}y_{1245} - y_{25}y_{1345} + y_{15}y_{2345}\\
 y_{15}y_{23} - y_{13}y_{25} + y_{12}y_{35} - y_{0}y_{1235} \\
 y_{14}y_{23} - y_{13}y_{24} + y_{12}y_{34} - y_{0}y_{1234} \\
 y_{1235}, y_{1234}
\end{align*}

\noindent
The resolution of $I_{20}$: 

\[R \leftarrow R^{9} \leftarrow R^{26} \leftarrow R^{37} \leftarrow R^{30} \leftarrow R^{14} \leftarrow R^{3} \leftarrow 0 \]

\noindent
 coincides with the tensor product of the Koszul complex on the variables $y_{1235}, y_{1234}$ and the  resolution of the variety of complexes studied by Miller in \cite{kustin1993} and defined by the ideal 
 
 \[I(X,Y|Y) = \operatorname{ideal}(X*Y)+\operatorname{minors}(2,Y)\]
 
 where 
 
\[
X = \begin{pmatrix}
y_{34} & y_{24} & y_{14} \\
y_{35} & y_{25} & y_{15}
\end{pmatrix}
\hbox{ and }
Y = \begin{pmatrix}
y_{12} & y_{1245} \\
-y_{13} & -y_{1345}\\
y_{23} & y_{2345}
\end{pmatrix}
\]

\noindent and where  $\operatorname{ideal}(X*Y)$ is the ideal generated by the entries of the matrix $X*Y$ and $\operatorname{minors}(2,Y)$ is the ideal of $2 \times 2$ minors of $Y$.

%
%
\begin{remark}
Since $g = 4, \max n_{gj} = 6, \min n_{1j} =2$,  Criterion (\ref{licci}) is satisfied. Hence, $I_{20}$ is not in the linkage class of a complete intersection.
\end{remark}

\subsection{($I_{19}$, $E_6$)}\label{i19e6} 
$I_{19}$ is a complete intersection of type $(1^4)$ in  the ideal generated by  $2\times 2$ minors of the $2\times 4$ matrix

$$N_{19} = \begin{pmatrix}{y}_{2345}&
       {y}_{23}&
       {y}_{24}&
       {y}_{25}\\
       {y}_{1345}&
       {y}_{13}&
       {y}_{14}&
       {y}_{15}\\
       \end{pmatrix}$$      
Tensoring the Eagon--Northcott complex of $N_{19}$, $R \leftarrow R^{6}(-2) \leftarrow R^{8}(-3) \leftarrow R^{3} (-4)\leftarrow 0$, with the Koszul complex on $4$ variables $y_{12}, y_{1245}, y_{1235}, y_{1234}$ we get the resolution of $I_{19}$:       
$$R \leftarrow R^{10} \leftarrow R^{38} \leftarrow R^{75} \leftarrow R^{85} \leftarrow R^{56} \leftarrow R^{20} \leftarrow R^{3} \leftarrow 0,$$
$P(19) =(1+3T)(1-T)^3 = -3T^{4} + 8T^{3} - 6T^{2}+1$. 

\subsection{($I_{18}$, $E_6$)}\label{i18e6} 
$I_{18}$ is resolved by Koszul complex on $8$ variables $y_{12}, y_{1245}, y_{1235}, y_{1234}, y_{1345}, y_{13}, y_{14}, y_{15}$, and hence is non-singular.

\subsection{($I_{17}$, $E_6$)}\label{i17e6} 
$I_{17}$ is generated by ten quadrics. The resolution is $$R \leftarrow R^{10}(-2) \leftarrow R^{16}(-3) \leftarrow R^{16}(-5) \leftarrow R^{10}(-6) \leftarrow R(-8) \leftarrow 0,$$

\noindent
which coincides with the expected resolution since \[P(I_{17}) = (1+5T+5T^2+T^3)(1-T)^5 = -T^8 +10T^6 -16T^5 +16T^3 -10T^2  +1.\]
This ideal defines the so-called variety of pure spinors (this very same description can be found in 
\begin{center}
https://math.galetto.org/orbits/index.html).
\end{center}

\noindent
It is the closure of the highest weight orbit in $V(\omega_{4},D_{5})$. The set of nodes $x$ in the Bruhat graph associated to $V(\omega_{1}, E6)$ which satisfy $x \leq x_{17}$,(without counting the identity, which we have set to be equal to one) are precisely the nodes corresponding to  $V(\omega_{4},D5)$ in the branching
 \[\operatorname{res}^{E6}_{D5}V(\omega_{1},E6) = V(0,D5)\oplus V(\omega_{1},D5) \oplus V(\omega_{5},D5).\]
 
\subsection{($I_{16}$, $E_6$)}\label{i16e6} 
$I_{16}$ is a a complete intersection of type $(1)$ in  $I_{17}$: we obtain its resolution 
\[R \leftarrow R^{11} \leftarrow R^{26} \leftarrow R^{32} \leftarrow R^{26} \leftarrow R^{11} \leftarrow R \leftarrow 0 \]
by tensoring the one for $I_{17}$ with the Koszul complex on $y_{1234}$.

\subsection{($I_{15}$, $E_6$)}\label{i15e6} 
$I_{15}$ is a complete intersection of type $(1^2)$ in  $I_{17}$; its resolution 
\[R \leftarrow R^{12} \leftarrow R^{37} \leftarrow R^{58} \leftarrow R^{58} \leftarrow R^{37} \leftarrow R^{12} \leftarrow R \leftarrow 0 \]
is obtained by tensoring the one for $I_{17}$ with the Koszul complex in the variables 
 $y_{1234}, y_{1235}$.   

\subsection{($I_{14}$, $E_6$)}\label{i14e6} 
$I_{14}$ is a complete intersection of type $(1^5)$ in  the ideal of $4\times 4$ Pfaffians of the skew-symmetric matrix 
\[ M_{14} = 
\begin{pmatrix}
0 & y_{45} & y_{35} & y_{25} &y_{15} \\
-y_{45} & 0 & y_{34} & y_{24} & y_{14}\\
-y_{35} & -y_{34} & 0 & y_{23} & y_{13} \\
-y_{25} & -y_{24} & -y_{23} & 0 & y_{12} \\
-y_{15}  & -y_{14} & -y_{13} & -y_{12} & 0
\end{pmatrix}.
\]
The resolution of $I_{14}$, 
\[ R \leftarrow R^{10} \leftarrow R^{40} \leftarrow R^{86} \leftarrow R^{110} \leftarrow R^{86} \leftarrow R^{40} \leftarrow R^{10} \leftarrow R \leftarrow 0, \]
is obtained by tensoring together the resolution 
$
R \leftarrow R^{5}(-2) \leftarrow R^{5}(-3) \leftarrow R(-5) \leftarrow 0
$
of the Pfaffian ideal with the Koszul complex in the five variables   $y_{1234}, y_{1235},
 y_{1245},y_{1345},y_{2345}$, $P(I_{14}) = (1+3T+T^2)(1-T)^3 = -T^{5} + 5T^3 - 5T^2 + 1$.  
 
\subsection{($I_{13}$, $E_6$)}\label{i13e6} 
$I_{13}$ 
Its resolution is  the tensor product of the Koszul complex in the variables $y_{1234}, y_{1235}, y_{1245},  y_{12}$ and the resolution

\[R \leftarrow R^{8} \leftarrow R^{12} \leftarrow R^{7} \leftarrow R^{2} \leftarrow 0 \]

with Betti table

\[\begin{matrix}
&0&1&2&3&4\\
\text{total:}&1&8&12&7&2\\
\text{0:}&1&\text{.}&\text{.}&\text{.}&\text{.}\\
\text{1:}&\text{.}&8&12&3&\text{.}\\
\text{2:}&\text{.}&\text{.}&\text{.}&4&2\\
     \end{matrix}
\]

 of the ideal 

 \[I(X,Y|Y) = \operatorname{ideal}(X*Y)+\operatorname{minors}(2,Y)\]
 
 where 
 
\[
X = \begin{pmatrix}
- y_{34} & y_{35} & -y_{45} & y_{0}
\end{pmatrix}
\hbox{ and }
Y = \begin{pmatrix}
-y_{25} & y_{15} \\
-y_{24} & y_{14} \\
-y_{23} & y_{13} \\
-y_{2345} & y_{1345}
\end{pmatrix}
\]

\noindent and where  $\operatorname{ideal}(X*Y)$ is the ideal generated by the entries of the matrix $X*Y$ and $\operatorname{minors}(2,Y)$ is the ideal of $2 \times 2$ minors of $Y$. We have

 \[P(13) = (1+4T + 2T^2)(1-T)^4 = 2T^6 -4T^5 -3T^4 + 12T^3 - 8T^2 + 1.\]

The ideal $I(X,Y|Y)$ is resolved by
\[
R \leftarrow R^8(-2) \leftarrow R^{12}(-3) \leftarrow R^{3}(-4) \oplus R^4(-5) \leftarrow R^2 (-6) \leftarrow 0. 
\] 

\begin{remark}
Since the codimension of $I_{13}$ is $4$, by Criterion \ref{licci} we have that since $g = 4, \max n_{gj} = 6, \min n_{1j} =2 $, (\ref{crit}) is satisfied. Hence, Criterion \ref{licci} implies that $(R/I_{13})_{(\underline{x})}$ is not in the linkage class of a complete intersection. 
\end{remark}

\subsection{($I_{12}$, $E_6$)}\label{j26e6} 
$I_{26}$ We have $P(I_{12}) = (1+3T+T^2)(1-T)^3 = -T^{5} + 5T^3 - 5T^2 + 1$. The expected resolution
\[ 
R \leftarrow R^{5}(-2) \leftarrow R^{5}(-3) \leftarrow R(-5) \leftarrow 0 
\]
is the resolution of the codimension 3 Gorenstein ideal generated by the quadrics $- y_{24}y_{13} + y_{23}y_{14},- y_{25}y_{13} + y_{15}y_{23},- y_{25}y_{14} + y_{24}y_{15}, y_{45}y_{13} - y_{35}y_{14} + y_{34}y_{15}, y_{45}y_{23} - y_{35}y_{24} + y_{34}y_{25}$, which are the Pfaffians of the sub-maximal minors of the matrix $M_{23}$ with the variable $y_{45}$ set to zero. The resolution of $I_{12}$ is obtained by tensoring the expected resolution with the Koszul complex in the variables 
$
y_{1234}, y_{1235}, y_{1245}, y_{1345}, y_{12}, y_{2345}.$

\subsection{($I_{11}$, $E_6$)}\label{i11e6} 
$I_{11}$ is given by a complete intersection of type $(1^{8},2)$. Its resolution is the Koszul complex in the elements $-y_0y_{2345} + y_{45}y_{23} - y_{35}y_{24} + y_{34}y_{25}, y_{1234},y_{1235},  y_{1245},y_{1345}, y_{12}, y_{13},y_{14},y_{15},$
$P(I_{11}) = (1+T)(1-T) = 1 - T^{2}$.

\subsection{($I_{10}$, $E_6$)}\label{i10e6} 
$I_{10}$ is given by a complete intersection of type $(1^{8})$ in the Eagon-Northcott ideal of $2 \times 2$ minors of the matrix
\[
 M_{10} = 
 \begin{pmatrix} 
 y_{34} & y_{24} & y_{14} \\
 y_{35} & y_{25} &y_{15}
 \end{pmatrix}.
\]
We have $P(I_{10}) = (1+2T)(1-T)^2 = 2T^3 - 3T^2 +1$. The resolution of $I_{10}$ is the tensor product of the Eagon-Northcott complex associated to $M_{10}$, $R \leftarrow R^3(-2) \leftarrow R^2(-3) \leftarrow 0$, and the Koszul complex in the variables 
$y_{23},y_{13},y_{12},y_{2345},y_{1345},y_{1245},y_{1235},y_{1234}.
$

\subsection{($I_{9}$, $E_6$)}\label{i9e6} 
$I_{9}$ is a complete intersection of type $(1^9,2)$.  The resolution of the ideal $I_{9}'$ generated by the quadric 
$y_{45}y_{23} - y_{35}y_{24} + y_{34}y_{25}$,  $R \leftarrow R(-2) \leftarrow 0$, coincides with the expected resolution, since we have 
\[P(I_{11}) = (1+T)(1-T) = 1-T^2.\] 
The resolution of $I_{9}$ is  the Koszul complex in the elements
$y_{45}y_{23} - y_{35}y_{24} + y_{34}y_{25}, y_{15},y_{14},y_{13},y_{12},y_{2345},y_{1345},y_{1245},y_{1235},y_{1234}.$  

\subsection{($I_{8}$, $E_6$)}\label{i8e6} $I_{8}$ is non-singular.

\subsection{($I_{7}$, $E_6$)}\label{i7e6} 
$I_{7}$ is a complete intersection of type $(1^{10})$ in  the ideal generated by the quadric $- y_{35}y_{24} + y_{34}y_{25}$, whose resolution coincides with the expected resolution since we have $P(I_7) = (1+T)(1-T) = 1-T^2$. The resolution of $I_{7}$
\[
R \leftarrow R^{11} \leftarrow R^{55} \leftarrow R^{165} \leftarrow R^{330} \leftarrow R^{462} \leftarrow R^{462} \leftarrow R^{330} \leftarrow R^{165} \leftarrow R^{55} \leftarrow R^{11} \leftarrow R^{11} \leftarrow R \leftarrow 0
\]
is the Koszul complex in the elements
$- y_{35}y_{24} + y_{34}y_{25}, y_{15},y_{14},y_{23},y_{13},y_{12},y_{2345},y_{1345},y_{1245},y_{1235},y_{1234}.$

\subsection{($I_{6}$--$I_{1}$, $E_6$)}\label{j61e6} Starting from $I_6$ all ideals have reduced Hilbert series equal to $1/(1-T)^{\dim}$, and thus are non-singular. The resolutions of these ideals $I$ are Koszul complexes in $\dim I$ variables.

%



\newpage
\newgeometry{margin=1cm}

\section{Ideals of Schubert varieties in $G(E_7)/P_7$}

\

$$
\begin{array}{l}
J_{55} = (Q) \\
J_{54} = (Q, f_{1}) \\
J_{53} = (Q, f_{1}, f_{2}) \\
J_{52} = (Q, f_{1}, f_{2},f_{3}) \\
J_{51} = (Q, f_{1}, f_{2},f_{3}, f_{4}, f_{6}, f_{7}) \\
J_{50} = (Q, f_{1}, f_{2},f_{3}, f_{4}, f_{5}) \\
J_{49} = (Q, f_{1}, f_{2},f_{3}, f_{4}, f_{5}, f_{6}, f_{8}, f_{10}, f_{12}, f_{15}, x_{27}) \\
J_{48} = (Q, f_{1}, f_{2},f_{3}, f_{4}, f_{5}, f_{6}, f_{7}) \\
J_{47} = (Q, f_{1}, f_{2},f_{3}, f_{4}, f_{5}, f_{6}, f_{7}, f_{8}, f_{10},f_{12},f_{15},x_{27}) \\
J_{46} = (Q, f_{1}, f_{2},f_{3}, f_{4}, f_{5}, f_{6}, f_{7}, f_{8}, f_{9}) \\
J_{45} = (Q, f_{1}, f_{2},f_{3}, f_{4}, f_{5}, f_{6}, f_{7}, f_{8}, f_{9}, f_{10}, f_{12}, f_{15}, x_{27}) \\
J_{44} = (Q, f_{1}, f_{2},f_{3}, f_{4}, f_{5}, f_{6}, f_{7}, f_{8}, f_{9}, f_{10}, f_{11}, f_{13}) \\
J_{43} = (Q, f_{1}, f_{2},f_{3}, f_{4}, f_{5}, f_{6}, f_{7}, f_{8}, f_{9}, f_{10}, f_{11}, f_{12}, f_{14}, f_{15}, f_{17}, x_{27}, x_{26}) \\
J_{42} = (Q, f_{1}, f_{2},f_{3}, f_{4}, f_{5}, f_{6}, f_{7}, f_{8}, f_{9}, f_{10}, f_{11}, f_{12}, f_{13}, f_{15}, x_{27}) \\
J_{41} = (Q, f_{1}, f_{2},f_{3}, f_{4}, f_{5}, f_{6}, f_{7}, f_{8}, f_{9}, f_{10}, f_{11}, f_{12}, f_{13}, f_{14}, f_{16}, f_{18}, f_{20}) \\
\end{array}
$$

$$
\begin{array}{l}
J_{40} = (Q, f_{1}, f_{2},f_{3}, f_{4}, f_{5}, f_{6}, f_{7}, f_{8}, f_{9}, f_{10}, f_{11}, f_{12}, f_{13}, f_{14}, f_{15}, f_{17}, x_{27}, x_{26}) \\
J_{39} = (Q, f_{1}, f_{2},f_{3}, f_{4}, f_{5}, f_{6}, f_{7}, f_{8}, f_{9}, f_{10}, f_{11}, f_{12}, f_{13}, f_{14}, f_{15}, f_{16}, f_{18}, f_{20}, x_{27}) \\
J_{38} = (Q, f_{1}, f_{2},f_{3}, f_{4}, f_{5}, f_{6}, f_{7}, f_{8}, f_{9}, f_{10}, f_{11}, f_{12}, f_{13}, f_{14}, f_{15}, f_{16}, f_{17}, f_{19}, x_{27}, x_{26}, x_{25}) \\
J_{37} = (Q, f_{1}, f_{2},f_{3}, f_{4}, f_{5}, f_{6}, f_{7}, f_{8}, f_{9}, f_{10}, f_{11}, f_{12}, f_{13}, f_{14}, f_{15}, f_{16}, f_{17}, f_{18}, f_{20}, x_{27},x_{26}) \\
J_{36} = (Q, f_{1}, f_{2},f_{3}, f_{4}, f_{5}, f_{6}, f_{7}, f_{8}, f_{9}, f_{10}, f_{11}, f_{12}, f_{13}, f_{14}, f_{15}, f_{16}, f_{17}, f_{18}, f_{19}, f_{21}, f_{23}, x_{27}, x_{26}, x_{25}, x_{24},x_{22},x_{20}) \\
J_{35} = (Q, f_{1}, f_{2},f_{3}, f_{4}, f_{5}, f_{6}, f_{7}, f_{8}, f_{9}, f_{10}, f_{11}, f_{12}, f_{13}, f_{14}, f_{15}, f_{16}, f_{17}, f_{18}, f_{19}, f_{20}, x_{27},x_{26},x_{25})\\
J_{34} = (Q, f_{1}, f_{2},f_{3}, f_{4}, f_{5}, f_{6}, f_{7}, f_{8}, f_{9}, f_{10}, f_{11}, f_{12}, f_{13}, f_{14}, f_{15}, f_{16}, f_{17}, f_{18}, f_{19}, f_{20}, f_{21}, f_{23}, x_{27}, x_{26}, x_{25}, x_{24}, x_{22}, x_{20})
\\
J_{33} = (Q, f_{1} \ra f_{22}, x_{27},x_{26},x_{25},x_{24},x_{23})
\end{array}$$

$$
\begin{array}{l}
J_{32} = (Q,f_{1}\ra f_{23}, x_{27},x_{26},x_{25},x_{24},x_{23},x_{22},x_{20}) \\ 

J_{31} = (Q,f_{1} \ra f_{24}, x_{27},x_{26},x_{25},x_{24},x_{23},x_{22},x_{21},x_{20},x_{18}) \\  

J_{30} = (Q,f_{1} \ra f_{25}, x_{27},x_{26},x_{25},x_{24},x_{23},x_{22},x_{21},x_{20},x_{18},x_{16},x_{13}) \\ 

J_{29} = (Q,f_{1} \ra f_{26}, x_{27},x_{26},x_{25},x_{24},x_{23},x_{22},x_{21},x_{20},x_{19},x_{18},x_{17},x_{16},x_{14},x_{13},x_{11},x_{9},x_{7}) \\  


J_{28} = (Q,f_{1} \ra f_{27}) \\ 

J_{27} = (Q,f_{1} \ra f_{27},x_{27}) \\ 

J_{26} = (Q,f_{1} \ra f_{27},x_{27},x_{26}) \\ 

J_{25} = (Q,f_{1} \ra f_{27},
x_{27},x_{26},x_{25}) \\ 

J_{24} = (Q,f_{1} \ra f_{27}, x_{27},x_{26},x_{25},x_{24},x_{22},x_{20}) \\ 

J_{23} = (Q,f_{1} \ra f_{27},
x_{27},x_{26},x_{25},x_{24},x_{23}) \\ 

J_{22} = (Q,f_{1} \ra f_{27}, x_{27},x_{26},x_{25},x_{24},x_{23},x_{22}) \\

J_{21} = (Q,f_{1} \ra f_{27}, x_{27},x_{26},x_{25},x_{24},x_{23},x_{22},x_{21},x_{19},x_{17},x_{15}) \\

J_{20} = (Q,f_{1} \ra f_{24}, x_{27},x_{26},x_{25},x_{24},x_{23},x_{22},x_{21},x_{20},x_{18}) \\

J_{19} = (Q,f_{1} \ra f_{27}, x_{27},x_{26},x_{25},x_{24},x_{23},x_{22},x_{21},x_{20},x_{19},x_{17},x_{15}) \\ 

J_{18} = (Q,f_{1} \ra f_{27}, x_{27},x_{26},x_{25},x_{24},x_{23},x_{22},x_{21},x_{20},x_{19},x_{18},x_{16},x_{13}) \\ 

J_{17} = (Q,f_{1} \ra f_{27}, 
x_{27},x_{26},x_{25},x_{24},x_{23},x_{22},x_{21},x_{20},x_{19},x_{18},x_{17},x_{15}) \\ 
\\
\end{array}
$$

$$
\begin{array}{l}
J_{16} =(Q,f_{1} \rightarrow f_{27},x_{27},x_{26},x_{25},x_{24},x_{23},x_{22},x_{21},x_{20},x_{19},x_{18},x_{17},x_{16},x_{14},x_{13},x_{11},x_{9},x_{7}) \\

J_{15} =(Q,f_{1} \rightarrow f_{27},x_{27},x_{26},x_{25},x_{24},x_{23},x_{22},x_{21},x_{20},x_{19},x_{18},x_{17},x_{16},x_{15},x_{13}) \\

J_{14} =(Q,f_{1} \rightarrow f_{27},x_{27},x_{26},x_{25},x_{24},x_{23},x_{22},x_{21},x_{20},x_{19},x_{18},x_{17},x_{16},x_{15},x_{14},x_{12}) \\

J_{13} =(Q,f_{1} \rightarrow f_{27},x_{27},x_{26},x_{25},x_{24},x_{23},x_{22},x_{21},x_{20},x_{19},x_{18},x_{17},x_{16},x_{15},x_{14},x_{13},x_{11},x_{9},x_{7}) \\

J_{12} =(Q,f_{1} \rightarrow f_{27},x_{27},x_{26},x_{25},x_{24},x_{23},x_{22},x_{21},x_{20},x_{19},x_{18},x_{17},x_{16},x_{15},x_{14},x_{13},x_{12}) \\

J_{11} =(Q,f_{1} \rightarrow f_{27},x_{27},x_{26},x_{25},x_{24},x_{23},x_{22},x_{21},x_{20},x_{19},x_{18},x_{17},x_{16},x_{15},x_{14},x_{13},x_{12},x_{11},x_{9},x_{7}) \\

J_{10} =(Q,f_{1} \rightarrow f_{27},x_{27},x_{26},x_{25},x_{24},x_{23},x_{22},x_{21},x_{20},x_{19},x_{18},x_{17},x_{16},x_{15},x_{14},x_{13},x_{12},x_{11},x_{10}) \\

J_{9} =(Q,f_{1} \rightarrow f_{27},x_{27},x_{26},x_{25},x_{24},x_{23},x_{22},x_{21},x_{20},x_{19},x_{18},x_{17},x_{16},x_{15},x_{14},x_{13},x_{12},x_{11},x_{10},x_{9},x_{7}) \\

J_{8} =(Q,f_{1} \rightarrow f_{27},x_{27},x_{26},x_{25},x_{24},x_{23},x_{22},x_{21},x_{20},x_{19},x_{18},x_{17},x_{16},x_{15},x_{14},x_{13},x_{12},x_{11},x_{10},x_{9},x_{8},x_{5}) \\

J_{7} =(Q,f_{1} \rightarrow f_{27},x_{27},x_{26},x_{25},x_{24},x_{23},x_{22},x_{21},x_{20},x_{19},x_{18},x_{17},x_{16},x_{15},x_{14},x_{13},x_{12},x_{11},x_{10},x_{9},x_{8},x_{7},x_{5}) \\

J_{6} =(Q,f_{1} \rightarrow f_{27},x_{27} \rightarrow x_{6}) \\

J_{5} =(Q,f_{1} \rightarrow f_{27},x_{27} \rightarrow x_{5}) \\

J_{4} =(Q,f_{1} \rightarrow f_{27},x_{27} \rightarrow x_{4}) \\

J_{3} =(Q,f_{1} \rightarrow f_{27},x_{27} \rightarrow x_{3}) \\

J_{2} =(Q,f_{1} \rightarrow f_{27},x_{27} \rightarrow x_{2}) \\

J_{1} =(Q,f_{1} \rightarrow f_{27},x_{27} \rightarrow x_{1}) \\
\end{array}
$$

\noindent
\subsubsection{Resolutions}
${\renewcommand{\arraystretch}{1.2}
\begin{array}{lllL}
\mbox{ideal} & \dim  & \mbox{numerator of HS} & \mbox{short description} \\
\hline
J_{55} & 26  &  -T^{3}+1 & c.i. of type $(3)$\\
\hline

J_{54} & 25 & 1+2 T+2 T^{2}+T^{3} & c.i. of type $(2,3)$  \\
\hline

J_{53} & 24 & 1+3 T+4 T^{2}+3 T^{3}+T^{4} & c.i. of type $(2,2,3)$  \\
\hline

J_{52} & 23 & 1+4 T+7 T^{2}+7 T^{3}+4 T^{4}+T^{5} & c.i. of type $(2,2,2,3)$ \\
\hline

J_{51} & 22 & 1+5 T+9 T^{2}+5 T^{3}+T^{4}  &  Huneke-Ulrich ideal of deviation 2\\
\hline

J_{50} & 22 & 1+5 T+10 T^{2}+9 T^{3}+2 T^{4} & Linked to $J_{51}$; codim. 5 complete intersection  \\
\hline

J_{49} & 21 & 1+5 T+5 T^{2}+T^{3}   &c.i. if type $(1)$ in variety  of pure spinors \\
\hline

J_{48} & 21 & 1+6 T+14 T^{2}+14 T^{3}+6 T^{4}+T^{5} & complete intersection of type $(2)$ in $J_{51}$  \\
\hline

J_{47} & 20 & 1+6T+10T^{2}+6T^{3} + T^{4} & c. i. of type $(1,2)$ in variety of pure spinors   \\
\hline

J_{46} & 20 & 1+7 T+19 T^{2}+22 T^{3}+10 T^{4}+T^{5} & ?  \\
\hline

J_{45} & 19 & 1+7 T+16 T^{2}+16 T^{3}+7 T^{4}+T^{5} &  c. i. of type $(2,2,1)$in variety of pure spinors  \\
\hline

J_{44} & 19 & 1+8 T+24 T^{2}+28 T^{3}+10 T^{4}+T^{5} & ?  \\
\hline

J_{43} & 18 &  1+7 T+13 T^{2}+7 T^{3}+T^{4} & ?  \\
\hline

J_{42} & 18 &  1+8 T+22 T^{2}+25 T^{3}+10 T^{4}+T^{5}  &  ? \\
\hline

J_{41} & 18 & 1+9 T+28 T^{2}+28 T^{3}+10 T^{4}+T^{5}  & ? \\
\hline
\hyperref[j39e7]{J_{39}} & 17 & 1+9 T+27 T^{2}+28 T^{3}+10 T^{4}+T^{5} & ?
 \\
\hline
\hyperref[j38e7]{J_{38}} & 16 & 1+8 T+18 T^{2}+13 T^{3}+2 T^{4} & ?
 \\
\hline
\hyperref[j37e7]{J_{37}} & 16 & 1+9 T+26 T^{2}+27 T^{3}+10 T^{4}+T^{5} & ?
\\
\hline
\hyperref[j36e7]{J_{36}} & 15 & 1+6 T+6 T^{2}+T^{3} 
& c.i. of type $(1^7)$ in ideal of $4 \times 4$ pfaff. in  $6\times 6$ skew-symm. \\
\hline
\hyperref[j35e7]{J_{35}} & 15 & 1+9 T+25 T^{2}+25 T^{3}+9 T^{4}+T^{5} & ?
 \\
\hline
\hyperref[j34e7]{J_{34}} & 14 & 1+8 T+16 T^{2}+8 T^{3}+T^{4} & ?
 \\
\hline
\hyperref[j33e7]{J_{33}} & 14 & 1+8 T+16 T^{2}+8 T^{3}+T^{4} & ?
\\
\hline
\hyperref[j32e7]{J_{32}} & 13 & 1+7 T+11 T^{2}+4 T^{3} & ?
 \\
\hline
\hyperref[j31e7]{J_{31}} & 12 & 1+6 T+6 T^{2} & ?
\\
\hline
\hyperref[j30e7]{J_{30}} & 10 & 1+4T & c.i. of type $(1^{12})$ in $2\times 2$ minors of a $2 \times 5$ generic matrix. \\
\hline
\hyperref[j29e7]{J_{29}} & 10 & 1 & non-singular   \\
\hline 
\hyperref[j28e7]{J_{28}} & 17 & 1+10 T+28 T^{2}+28 T^{3}+10 T^{4}+T^{5} 
& $(E_{6},P_{1})$\\
\hline 

\hyperref[j27e7]{J_{27}} & 16 & 1+10 T+28 T^{2}+28 T^{3}+10 T^{4}+T^{5} 
& c.i. of type $(1)$ in $(E_{6},P_{1})$\\
\hline 

\hyperref[j26e7]{J_{26}} & 15 & 1+10 T+28 T^{2}+28 T^{3}+10 T^{4}+T^{5} 
& c.i. of type $(1,1)$ in  $(E_{6},P_{1})$\\
\hline 

\hyperref[j25e7]{J_{25}} & 14 & 1+10 T+28 T^{2}+28 T^{3}+10 T^{4}+T^{5} 
& c.i. in $(E_{6},P_{1})$ \\
\hline 

\hyperref[j24e7]{J_{24}} & 13 & 1+8 T+15 T^{2}+8 T^{3}+T^{4} 
& ?\\
\hline 

\hyperref[j23e7]{J_{23}} & 13 & 1+9 T+20 T^{2}+13 T^{3}+2 T^{4} 
& ?\\
\hline 

\hyperref[j22e7]{J_{22}} & 12 & 1+9 T+20 T^{2}+13 T^{3}+2 T^{4} 
& ? \\
\hline 

\hyperref[j21e7]{J_{21}} & 12 & 1+5 T+5 T^{2}+T^{3} 
& $V(\omega_{5},D_{5})$ \\
\hline 

\hyperref[j20e7]{J_{20}} & 11 & 1+6T+6 T^{2} 
& ? \\
\hline 

\hyperref[j19e7]{J_{19}} & 11 & 1+5 T+5 T^{2}+T^{3} 
& $V(\omega_{5},D_{5})$ \\
\hline 

\hyperref[j18e7]{J_{18}} & 10 & 1+5 T+3 T^{2} 
& c.i. of type $(1^{12})$ in variety of complexes \\
\hline 

\hyperref[j17e7]{J_{17}} & 10 & 1+5 T+5 T^{2}+T^{3} 
& $V(\omega_{5},D_{5})$ \\
\hline 

\end{array}}
$

$${\renewcommand{\arraystretch}{1.2}
\begin{array}{lllL}
\mbox{ideal} & \dim & \mbox{numerator} & \mbox{short description}\\
\hline 
\hyperref[j16e7]{J_{16}} & 9 & {T+1} 
& c.i. of type $(1^{17},2)$ \\
\hline 
\hyperref[j15e7]{J_{15}} & 9 & {2T^2+4T+1} 
& c.i. of type $(1^{14})$ in the ideal associated to a matrix and a vector  \cite{kustin1993}\\
\hline 
\hyperref[j14e7]{J_{14}} & 9 & {T^2+3T+1} 
& c.i. of type $(1^{15})$ in Pf. of $5\times 5$ skew\\
\hline 
\hyperref[j13e7]{J_{13}} & 8 & {T+1} 
&c.i. of type $(1^{18},2)$ \\
\hline 
\hyperref[j12e7]{J_{12}} & 8 & {T^2+3T+1} 
& c.i. of type $(1^{16})$ in Pf. of $5\times 5$ skew \\
\hline 
\hyperref[j11e7]{J_{11}} & 7 & {T+1} 
&  c.i. of type $(1^{19},2)$ \\
\hline 
\hyperref[j10e7]{J_{10}} & 7 & {2T+1} 
& $2 \times 2$ minors of  a generic $2 \times 3$ matrix \\
\hline 
\hyperref[j9e7]{J_{9}} & 6 & {T+1} 
&  c.i. of type $(1^{20},2)$ \\
\hline 
\hyperref[j81e7]{J_{8}} & 6 & {1} & non-singular \\
\hline 
\hyperref[j81e7]{J_{7}} & 5 & {1}  & non-singular \\
\hline 
\hyperref[j81e7]{J_{6}} & 5 & {1}  & non-singular \\
\hline 
\hyperref[j81e7]{J_{5}} & 4 & {1}  & non-singular \\
\hline 
\hyperref[j81e7]{J_{4}} & 3 & {1} & non-singular \\
\hline 
\hyperref[j81e7]{J_{3}} & 2 & {1} & non-singular \\
\hline 
\hyperref[j81e7]{J_{2}} & 1 & {1} & non-singular \\
\hline 
\hyperref[j81e7]{J_{1}} & 0 & {1} & non-singular \\
\hline
\end{array}}
$$

\

\subsection{Descriptions}

For each ideal $J_{s}$ we will denote, whenever appropriate, by $J'_{s}$ the ideal such that $J_{s}$ is generated by $J'_{s}$ and coordinate hyperplane sections. 

\subsection{($J_{55}$, $E_7$)}\label{j55e7} 
$J_{55}= (Q)$ is resolved by Koszul complex on the cubic $Q$, $R \leftarrow R(-3) \leftarrow 0$ and  $P(J_{55}) = 1-T^{3}$.

\subsection{($J_{54}$, $E_7$)}\label{j54e7} 
$J_{54} = (Q, f_{1})$ is resolved by Koszul complex on the cubic $Q$ and the quadric $f_{1}$, $R \leftarrow R(-2)\oplus R(-3) \leftarrow R(-5) \leftarrow 0$; we have $P(J_{54}) = T^{5}-T^{3}-T^{2}+1$.

\subsection{($J_{53}$, $E_7$)}\label{j53e7} 
$J_{53} = (Q, f_{1}, f_{2})$ is a complete intersection resolved by Koszul complex on the cubic $Q$ and the two quadrics $f_{1}$ and $f_{2}$:
$$R \leftarrow R^{2}(-2) \oplus R(-3) \leftarrow R^{2}(-5) \oplus R(-4)  \leftarrow R(-7) \leftarrow 0;$$ 
$we have P(J_{53}) = -T^{7}+2 T^{5}+T^{4}-T^{3}-2 T^{2}+1$.

\subsection{($J_{52}$, $E_7$)}\label{j52e7} 
$J_{52} = (Q, f_{1}, f_{2},f_{3})$ is a complete intersection on the four generators $Q,f_1,f_2,f_3$ resolved by Koszul complex on the cubic $Q$ and the three quadrics $f_{25}$, $f_{26}$ and $f_{27}$:
$$R \la R^3(-2) \oplus R(-3) \la R^{3}(-4) \oplus R^3(-5) \la R(-6) \oplus R^3(-7) \la R(-9) \leftarrow 0; $$
$we have P(J_{52}) = T^{9}-3 T^{7}-T^{6}+3 T^{5}+3 T^{4}-T^{3}-3 T^{2}+1$. 

\subsection{($J_{51}$, $E_7$)}\label{j51e7} 
We have $J_{51} = (Q, f_{1}, f_{2},f_{3}, f_{4}, f_{6}, f_{7})$ and
  \[P(J_{51}) = (T^4+5T^3+9T^2+5T+1)(1-T)^5 = -T^9+6T^7-21T^5+21T^4-6T^2+1.\] 
  
  If we choose $Y = (x_{20},x_{22},x_{24},x_{25},x_{26},x_{27})$,
\[
X = \begin{pmatrix}
0 & {x}_{23} & -{x}_{21} & {x}_{19}  & -{x}_{17} & {x}_{15} \\
- {x}_{23} &  0 &  {x}_{18} & -{x}_{16} & {x}_{14} & -{x}_{12}\\
{x}_{21} &  -{x}_{18} &  0 &  {x}_{13} & -{x}_{11} & {x}_{10}\\
-{x}_{19} &  {x}_{16} &  -{x}_{13} &  0 &  {x}_{9} & -{x}_{8}\\
{x}_{17} & -{x}_{14} & {x}_{11} & -{x}_{9}  & 0 & {x}_{5}\\
-{x}_{15} & {x}_{12} & -{x}_{10} & {x}_{8} & -{x}_{5} & 0
\end{pmatrix},
\]
then $J_{51} = I_1(Y\cdot X) + \Pf(X)$, hence $J_{51}$ is a $7$-generated Huneke--Ulrich ideal of deviation two  studied in \cite[Example 2.7]{Kustin1986} with resolution

\[\mathbb{F}_{51} = R \leftarrow R^{7} \leftarrow R^{22} \leftarrow R^{22} \leftarrow R^{7} \leftarrow R \leftarrow 0 \]

of Betti table:
\[
\begin{matrix}
       &0&1&2&3&4&5\\
  \text{total:}&1&7&22&22&7&1\\
  \text{0:}&1&\text{.}&\text{.}&\text{.}&\text{.}&\text{.}\\
  \text{1:}&\text{.}&6&1&\text{.}&\text{.}&\text{.}\\
  \text{2:}&\text{.}&1&21&21&1&\text{.}\\
  \text{3:}&\text{.}&\text{.}&\text{.}&1&6&\text{.}\\
  \text{4:}&\text{.}&\text{.}&\text{.}&\text{.}&\text{.}&1
 \end{matrix}
\]
\

\subsection{($J_{50}$, $E_7$)}\label{j50e7} 
$J_{50} = (Q, f_{1}, f_{2},f_{3}, f_{4}, f_{5})$ is linked to $J_{51}$ by the five elements $Q,f_{1},f_{2},f_{3},f_{4}$.  \\
We have  \[P(J_{50}) = 2T^4+9T^3+10T^2+5T+1)(1-T)^5=-2T^9+T^8+15T^7-25T^6+4T^5+12T^4-T^3-5T^2+1.\]

\

\subsection{($J_{49}$, $E_7$)}\label{j49e7} 
$J_{49} = (Q, f_{1}, f_{2},f_{3}, f_{4}, f_{5}, f_{6}, f_{8}, f_{10}, f_{12}, f_{15}, x_{27})$ is a complete intersection of type $(1)$ in  in the variety of pure spinors, hence it is isomorphic to $I_{16}$.  \\
We have  $P(J_{49}) = (T^3+5T^2+5T+1)(1-T)^5 = -T^8+10T^6-16T^5+16T^3-10T^2+1$ and its resolution is

\[R \leftarrow R^{11} \leftarrow R^{26} \leftarrow R^{32} \leftarrow R^{26} \leftarrow R^{11} \leftarrow R \leftarrow 0, \]

\noindent
isomorphic to the resolution of $I_{16}$, i.e. the tensor product of the resolution of the variety of pure spinors with the Koszul complex in the one variable $x_{27}$. 
\

\subsection{($J_{48}$, $E_7$)}\label{j48e7} 
$J_{48} = (Q, f_{1}, f_{2},f_{3}, f_{4}, f_{5}, f_{6}, f_{7})
$\\

note $(T^5+6T^4+14T^3+14T^2+6T+1)(1-T)^6= T^{11}-7T^9+27T^7-21T^6-21T^5+27T^4-7T^2+1$.

\subsection{($J_{47}$, $E_7$)}\label{j47e7} 
$J_{47} = (Q, f_{1}, f_{2},f_{3}, f_{4}, f_{5}, f_{6}, f_{7}, f_{8}, f_{10},f_{12},f_{15},x_{27}$)
Is a complete intersection in the variety of pure spinors. Its resolution is the tensor product of $\mathbf{F}_{ps}$, the resolution of the variety of pure spinors, and the Koszul complex on the elements $x_{27} $ and $ {x}_{15}{x}_{20}-{x}_{12}{x}_{22}+{x}_{10}{x}_{24}-{x}_{8}{x}_{25}+{x}_{5}{x}_{26}$.

%

\subsection{($J_{46}$, $E_7$)}\label{j46e7} 
$J_{46} = (Q, f_{1}, f_{2},f_{3}, f_{4}, f_{5}, f_{6}, f_{7}, f_{8}, f_{9}$)

\begin{align*}
{x}_{22}{x}_{23}-{x}_{21}{x}_{24}+{x}_{19}{x}_{25}-{x}_{17}{x}_{26}+{x}_{15}{x}_{27},\\
{x}_{20}{x}_{23}-{x}_{18}{x}_{24}+{x}_{16}{x}_{25}-{x}_{14}{x}_{26}+{x}_{12}{x}_{27},\\
{x}_{20}{x}_{21}-{x}_{18}{x}_{22}+{x}_{13}{x}_{25}-{x}_{11}{x}_{26}+{x}_{10}{x}_{27},\\
{x}_{19}{x}_{20}-{x}_{16}{x}_{22}+{x}_{13}{x}_{24}-{x}_{9}{x}_{26}+{x}_{8}{x}_{27},\\
{x}_{17}{x}_{20}-{x}_{14}{x}_{22}+{x}_{11}{x}_{24}-{x}_{9}{x}_{25}+{x}_{5}{x}_{27},\\
{x}_{15}{x}_{20}-{x }_{12}{x}_{22}+{x}_{10}{x}_{24}-{x}_{8}{x}_{25}+{x}_{5}{x}_{26},\\
{x}_{18}{x}_{19}-{x}_{16}{x}_{21}+{x}_{13}{x}_{23}-{x}_{7}{x}_{26}+{x}_{6}{x}_{27},\\
{x}_{17}{x}_{18}-{x}_{14}{x}_{21}+{x}_{11}{x}_{23}-{x}_{7}{x}_{25}+{x}_{4}{x}_{27},\\
{x}_{15}{x}_{18}-{x}_{12}{x}_{21}+{x}_{10}{x}_{23}-{x}_{6}{x}_{25}+{x}_{4}{x}_{26},\\
{x}_{13}{x}_{14}{x}_{15}-{x}_{11}{x}_{15}{x}_{16}-{x}_{12}{x}_{13}{x}_{17}\\
+{x}_{10}{x}_{16}{x}_{17}+{x}_{11}{x}_{12}{x}_{19}-{x}_{10}{x}_{14}{x}_{19}\\
-{x}_{7}{x}_{8}{x}_{25}+{x}_{6}{x}_{9}{x}_{25}+{x}_{5}{x}_{7}{x}_{26}\\
-{x}_{4}{x}_{9}{x}_{26}-{x}_{5}{x}_{6}{x}_{27}+{x}_{4}{x}_{8}{x}_{ 27}
\end{align*}

\subsection{($J_{45}$, $E_7$)}\label{j45e7} 
$J_{45} = (Q, f_{1}, f_{2},f_{3}, f_{4}, f_{5}, f_{6}, f_{7}, f_{8}, f_{9}, f_{10}, f_{12}, f_{15}, x_{27}$)
is a complete intersection of type $(2,2)$ in the variety of pure spinors.

\subsection{($J_{44}$, $E_7$)}\label{j44e7} 
$J_{44} = (Q, f_{1}, f_{2},f_{3}, f_{4}, f_{5}, f_{6}, f_{7}, f_{8}, f_{9}, f_{10}, f_{11}, f_{13}$)

\begin{align*}
{x}_{22}{x}_{23}-{x}_{21}{x}_{24}+{x}_{19}{x}_{25}-{x}_{17}{x}_{26}+{x}_{15}{x}_{27},\\
{x}_{20}{x}_{23}-{x}_{18}{x}_{24}+{x}_{16}{x}_{25}-{x}_{14}{x}_{26}+{x}_{12}{x}_{27},\\
{x}_{20}{x}_{21}-{x}_{18}{x}_{22}+{x}_{13}{x}_{25}-{x}_{11}{x}_{26}+{x}_{10}{x}_{27},\\
{x}_{19}{x}_{20}-{x}_{16}{x}_{22}+{x}_{13}{x}_{24}-{x}_{9}{x}_{26}+{x}_{8}{x}_{27},\\
{x}_{17}{x}_{20}-{x}_{14}{x}_{22}+{x}_{11}{x}_{24}-{x}_{9}{x}_{25}+{x}_{5}{x}_{27},\\
{x}_{15}{x}_{20}-{x}_{12}{x}_{22}+{x}_{10}{x}_{24}-{x}_{8}{x}_{25}+{x}_{5}{x}_{26},\\
{x}_{18}{x}_{19}-{x}_{16}{x}_{21}+{x}_{13}{x}_{23}-{x}_{7}{x}_{26}+{x}_{6}{x}_{27},\\
{x}_{17}{x}_{18}-{x}_{14}{x}_{21}+{x}_{11}{x}_{23}-{x}_{7}{x}_{25}+{x}_{4}{x}_{27},\\
{x}_{15}{x}_{18}-{x}_{12}{x}_{21}+{x}_{10}{x}_{23}-{x}_{6}{x}_{25}+{x}_{4}{x}_{26},\\
{x}_{16}{x}_{17}-{x}_{14}{x}_{19}+{x}_{9}{x}_{23}-{x}_{7}{x}_{24}+{x}_{3}{x}_{27},\\
{x}_{15}{x}_{16}-{x}_{12}{x}_{19}+{x}_{8}{x}_{23}-{x}_{6}{x}_{24}+{x}_{3}{x}_{26},\\
{x}_{14}{x}_{15}-{x}_{12}{x}_{17}+{x}_{5}{x}_{23}-{x}_{4}{x}_{24}+{x}_{3}{x}_{25},\\
{x}_{9}{x}_{10}{x}_{23}-{x}_{8}{x}_{11}{x}_{23}+{x}_{5}{x}_{13}{x}_{23}\\
-{x}_{7}{x}_{10}{x}_{24}+{x}_{6}{x}_{11}{x}_{24}-{x}_{4}{x}_{13}{x}_{24}\\
+{x}_{7}{x}_{8}{x}_{25}-{x}_{6}{x}_{9}{x}_{25}+{x}_{3}{x}_{13}{x}_{25}\\
-{x}_{5}{x}_{7}{x}_{26}+{x}_{4}{x}_{9}{x}_{26}-{x}_{3}{x}_{11}{x}_{26}\\
+{x}_{5}{x}_{6}{x}_{27}-{x}_{4}{x}_{8}{x}_{27}+{x}_{3}{x}_{10}{x}_{27}
\end{align*}

%
%
%
%
%
\subsection{($J_{36}$, $E_7$)}\label{j36e7} 
This ideal is a complete intersection of type $(1^7)$ in  ideal generated by the  $4\times 4$ pfaffians of the skew-symmetric matrix below: 


\[
X = \begin{pmatrix}
0& x_{23} & -x_{21} & x_{19} & -x_{17}& x_{15} \\
- x_{23} &  0 &  x_{18} & -x_{16} & x_{14} & -x_{12} \\
x_{21} &  -x_{18} &  0 &  x_{13} & -x_{11} & x_{10} \\
-x_{19} & x_{16} &  -x_{13} &  0 &  x_{9} & -x_{8} \\
x_{17} & -x_{14} & x_{11} & -x_{9} & 0 & x_{5}\\
-x_{15} & x_{12} & -x_{10} & x_{8} & -x_{5} & 0
\end{pmatrix}
\]

%
%
%
%
%
\subsection{($J_{30}$, $E_7$)}\label{j30e7} 
$J_{30}$ is a complete intersection of type $(1^{12})$ in  the ideal $J_{30}'$ generated by all the $2 \times 2$ minors of the matrix
 
 \[X = 
 \begin{pmatrix}
x_{7} & -x_{9} & x_{11} & -x_{14} & x_{17} \\
-x_{6} & x_{8} & -x_{10} & x_{12} & -x_{15}
 \end{pmatrix}
 \]
 
 The resolution of $J_{30}$ is the tensor product $\mathbb{F}'_{30} \oplus \mathbb{K}_{12}$ where \
 
 \[\mathbb{F}'_{30} = R \leftarrow R^{10} \leftarrow R^{20} \leftarrow R^{15} \leftarrow R^{4} \leftarrow 0  \]
 
 is the resolution of $J_{30}'$ and $\mathbb{K}_{12}$ is the Koszul complex on the variables $x_{13},x_{16},x_{18},x_{19},x_{20},x_{21},x_{22},x_{23},x_{24},x_{25},x_{26},x_{27}$. 
%
%
%
%
%
%
%
%
%
%
%
\subsection{($J_{18}$, $E_7$)}\label{j18e7} 
$J_{18}$ is 
a complete intersection of type $(1^{12})$ in  in the variety of complexes

 \[I(X,Y|Y) = \operatorname{ideal}(X*Y)+\operatorname{minors}(2,Y)\]
 
 defined by $J_{18}'$, where 
 
 \[ X =
 \begin{pmatrix}
-x_{5} & -x_{4} & -x_{3} &-x_{2} & -x_{1}
 \end{pmatrix}
 \]

and

 \[ Y =
  \begin{pmatrix}
x_{7} & -x_{6}\\
-x_{9} & x_{8} \\
x_{11} & -x_{10} \\
-x_{14} & x_{12} \\
x_{17} & -x_{15}
 \end{pmatrix}
 \]

Its resolution is the tensor product  $\mathbb{F}'_{18} \oplus \mathbb{K}_{12}$ where \
 
 \[\mathbb{F}'_{18} = R \leftarrow R^{12} \leftarrow R^{25} \leftarrow R^{25} \leftarrow R^{14} \leftarrow R^3 \leftarrow 0  \]
 
 is the resolution of $J_{18}'$ and $\mathbb{K}_{12}$ is the Koszul complex on the variables $x_{13},x_{16},x_{18},x_{19},x_{20},x_{21},x_{22},x_{23},x_{24},x_{25},x_{26},x_{27}$. 


\subsection{($J_{16}$, $E_7$)}\label{j16e7} 
$J_{16}'
= \langle q_{16} \rangle 
= \langle x_{5}x_{6}-x_{4}x_{8}+x_{3}x_{10}-x_{2}x_{12}+x_{1}x_{15} \rangle$, 
$P(T) = {-T^2+1}$, and 
$\mathbb{F}'_{16} \, : \, R \leftarrow R(-2) \leftarrow 0$. \\
$J_{16}$ is resolved by the Koszul complex onthe 18 elements 
$$q_{16}, x_{27},x_{26},x_{25},x_{24},x_{23},x_{22},x_{21},x_{20},x_{19},x_{18},x_{17},x_{16},x_{14},x_{13},x_{11},x_{9},x_{7}.$$
\subsection{($J_{15}$, $E_7$)}\label{j15e7} 
 If we take the following $2 \times 4$ matrix and $4$-vector,
$$
X_{15} =\begin{pmatrix}{x}_{6}&
      {x}_{8}&
      {x}_{10}&
      {x}_{12}\\
      {x}_{7}&
      {x}_{9}&
      {x}_{11}&
      {x}_{14}\\
\end{pmatrix} \quad
\mbox{and} \quad
v_{15} = \begin{pmatrix} {-x_5},{x_4},{-x_3},{x_2} \end{pmatrix},$$ 
 $J_{15}'$ is the ideal associated to a matrix and a vector
 
  \[J_{15}' = \operatorname{ideal}(X_{15}\cdot v_{15})+\operatorname{minors}(2,X_{15})\]
 
described in \cite{kustin1993}. Note that $P(T)= {2T^6-4T^5-3T^4+12T^3-8T^2+1}$ and $J'_{15}$ is resolved by 
$$\mathbb{F}'_{15}=  R \leftarrow R^{8}(-2) \leftarrow R^{12}(-3) \leftarrow R^{3}(-4) \oplus R^4(-5) \leftarrow R^{2}(-6) \leftarrow 0.$$\\

In this case $J_{15}$ is a complete intersection of type $(1^{14})$ in  $J'_{15}$ and its resolution is therefore obtained by tensoring $\mathbb{F}_{15}$ with the Koszul complex on the 14 variables 
$$x_{27},x_{26},x_{25},x_{24},x_{23},x_{22},x_{21},x_{20},x_{19},x_{18},x_{17},x_{16},x_{15},x_{13}.$$

\subsection{($J_{14}$, $E_7$)}\label{j14e7} 
$J_{14}'
= \langle q_{14}^{1}, \ldots, q_{14}^{5} \rangle 
= \langle x_{9}x_{10}-x_{8}x_{11}+x_{5}x_{13}, -x_{7}x_{10}+x_{6}x_{11}-x_{4}x_{13}, x_{7}x_{8}-x_{6}x_{9}+x_{3}x_{13}, -x_{5}x_{7}+x_{4}x_{9}-x_{3}x_{11}, x_{5}x_{6}-x_{4}x_{8}+x_{3}x_{10}
\rangle$.
 We have $P(T)= {-T^5+5T^3-5T^2+1}$ and ${\mathbb{F}}'_{14} \, : \, $ {$R \leftarrow R^{5}(-2) \leftarrow R^{5}(-3) \leftarrow R(-5) \leftarrow 0$}.\\
$J_{14}'$ is generated by maximal Pfaffians of the $5\times 5$ (skewsymmetric) matrix 
$$
X_{14}= \begin{pmatrix}0&
      {x}_{3}&
      {x}_{4}&
      {x}_{6}&
      {x}_{7}\\
      {-{x}_{3}}&
      0&
      {x}_{5}&
      {x}_{8}&
      {x}_{9}\\
      {-{x}_{4}}&
      {-{x}_{5}}&
      0&
      {x}_{10}&
      {x}_{11}\\
      {-{x}_{6}}&
      {-{x}_{8}}&
      {-{x}_{10}}&
      0&
      {x}_{13}\\
      {-{x}_{7}}&
      {-{x}_{9}}&
      {-{x}_{11}}&
      {-{x}_{13}}&
      0\\
      \end{pmatrix}
      $$
Hence $J_{14}$ is a complete intersection of type $(1^{14})$ in  this Pfaffian ideal, and its resolution $\mathbb{F}_{14} = \mathbb{F}'_{14} \otimes \mathbb{K}_{15}$ is the tensor product of $\mathbb{F}_{14}'$ withe the Koszul complex $\mathbb{K}_{15}$ on the  variables
$$x_{27}, x_{26},x_{25},x_{24},x_{23},x_{22},x_{21},x_{20},x_{19},x_{18},x_{17},x_{16},x_{15},x_{14},x_{12}.$$ 

\subsection{($J_{13}$, $E_7$)}\label{j13e7} 
$J_{13}'= \langle q_{13} \rangle = \langle x_{5}x_{6}-x_{4}x_{8}+x_{3}x_{10}-x_{2}x_{12} \rangle$, $P(T)= {-T^2+1}$ and ${\mathbb{F}}'_{13} \, : \, $ $R \leftarrow R(-2) \leftarrow 0$.\\
$J_{13}$ is resolved by Koszul complex on the quadric $q_{13}$ tensored by Koszul complex on the 18 variables $$x_{27},x_{26},x_{25},x_{24},x_{23},x_{22},x_{21},x_{20},x_{19},x_{18},x_{17},x_{16},x_{15},x_{14},x_{13},x_{11},x_{9},x_{7}.$$

\subsection{($J_{12}$, $E_7$)}\label{j12e7} 
$J_{12}'
= \langle q_{12}^{1}, \ldots, q_{12}^{5} \rangle 
= \langle -{x}_{7} {x}_{8}+{x}_{6} {x}_{9},{x}_{5} {x}_{6}-{x}_{4}
      {x}_{8}+{x}_{3} {x}_{10},-{x}_{9} {x}_{10}+{x}_{8} {x}_{11},-{x}_{7}
      {x}_{10}+{x}_{6} {x}_{11},{x}_{5} {x}_{7}-{x}_{4} {x}_{9}+{x}_{3}
      {x}_{11} \rangle$, 
$P(T)= {-T^5+5T^3-5T^2+1}$ and 
${\mathbb{F}}'_{12} \, : \, $ {$R \leftarrow R^{5}(-2) \leftarrow R^{5}(-3) \leftarrow R(-5) \leftarrow 0$}.\\
$J_{12}'$ is a complete intersection of type $(1)$ in $J_{14}'$, in particular it is a complete intersection of type $1^{16}$ in the ideal of maximal Pfaffians of $X_{14}$. 
$J_{12}$ is resolved by tensoring it with the Koszul complex on 16 variables $$x_{27},x_{26},x_{25},x_{24},x_{23},x_{22},x_{21},x_{20},x_{19},x_{18},x_{17},x_{16},x_{15},x_{14},x_{13},x_{12}.$$

\subsection{($J_{11}$, $E_7$)}\label{j11e7} 
$J_{11}'
= \langle q_{11} \rangle 
= \langle x_{5}x_{6}-x_{4}x_{8}+x_{3}x_{10} \rangle$, 
$P(T)= {-T^2+1}$ and 
${\mathbb{F}}'_{11} \, : \, $ $R \leftarrow R(-2) \leftarrow 0$.\\
$J_{11}$ is resolved by Koszul complex on the quadric $q_{11}$ tensored by Koszul complex on 19 variables 
$$
x_{27},x_{26},x_{25},x_{24},x_{23},x_{22},x_{21},x_{20},x_{19},x_{18},x_{17},x_{16},x_{15},x_{14},x_{13},x_{12},x_{11},x_{9},x_{7}.
$$

\subsection{($J_{10}$, $E_7$)}\label{j10e7} 
$J_{10}'= \langle q_{10}^1,q_{10}^2,q_{10}^3 \rangle = \langle x_{7}x_{8}-x_{6}x_{9}, -x_{5}x_{7}+x_{4}x_{9}, x_{5}x_{6}-x_{4}x_{8} \rangle$, $P(T)= {2T^3-3T^2+1}$ and ${\mathbb{F}}'_{10} \, : \, $ $R \leftarrow R^{3}(-2) \leftarrow R^{2}(-3) \leftarrow 0$.\\
$J_{10}$ is resolved by Koszul complex on three quadrics $q_{10}^1,q_{10}^2,q_{10}^3$ tensored by Koszul complex on 18 variables $$x_{27},x_{26},x_{25},x_{24},x_{23},x_{22},x_{21},x_{20},x_{19},x_{18},x_{17},x_{16},x_{15},x_{14},x_{13},x_{12},x_{11},x_{10}.$$

\subsection{($J_{9}$, $E_7$)}\label{j9e7} 
$J_{9}'= \langle q_9 \rangle = \langle x_{5}x_{6}-x_{4}x_{8} \rangle$, $P(T)= {-T^2+1}$ and ${\mathbb{F}}'_{9} \, : \, $ $R \leftarrow R(-2) \leftarrow 0$.\\
$J_{16}$ is resolved by the Koszul complex on the  21 elements
$$x_{5}x_{6}-x_{4}x_{8},x_{27},x_{26},x_{25},x_{24},x_{23},x_{22},x_{21},x_{20},x_{19},x_{18},x_{17},x_{16},x_{15},x_{14},x_{13},x_{12},x_{11},x_{10},x_{9},x_{7}.$$

\subsection{($J_{8}$-$J_{1}$, $E_7$)}\label{j81e7} The ideals from $J_8$ to $J_1$ are non-singular.

\appendix
\label{sec:appendix}

\newgeometry{margin=1cm}
\section{Bruhat graphs and labellings}
\label{app:bruhat}
\FloatBarrier
 
\begin{figure}
\label{graphe7p7}
\includegraphics[scale = 0.6]{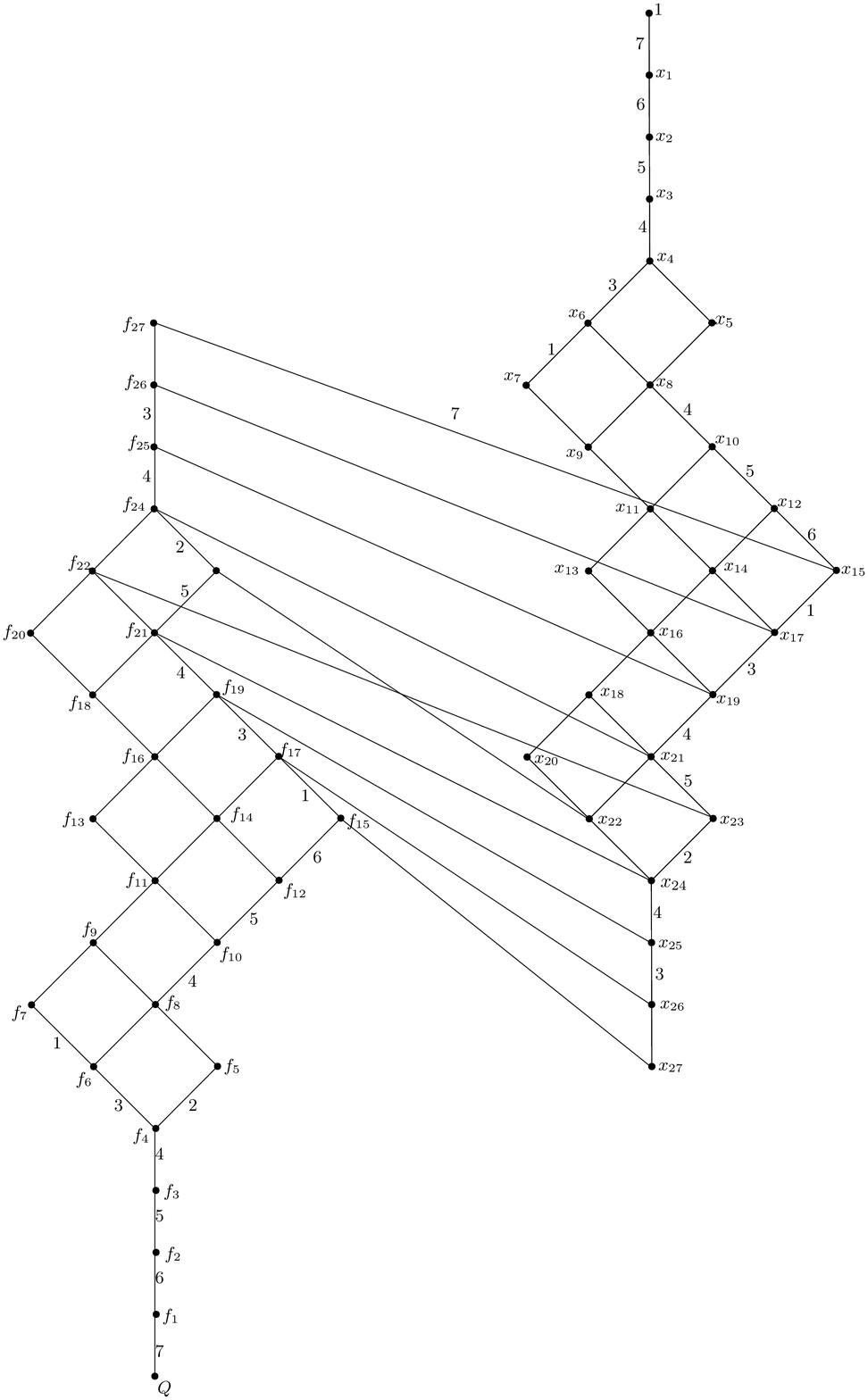}
\caption{Bruhat interval $W/W_{P_{7}}$ for $E_{7}$}
\end{figure}

\section{Weights and Symmetries}
\subsection{$E_{6}$}
$$
{\renewcommand{\arraystretch}{0.8}
\begin{array}{llllll}
\mbox{No}
&
\mbox{Weights of $V(\varpi_1)$}
&
\hbox{elt. in }
W/W_{P}
&
\mbox{Group of symmetries}
&
\mbox{Gorenstein}
&
dim
\\ 
\hline
\\
1
&
\begin{pmatrix}
1 & 0 & 0 & 0& 0\\
 & & 0 & &  
\end{pmatrix} 
& 
e
&
E_6
&
\hbox{yes (non-singular)}
&
0
\\
\\
2
&
\begin{pmatrix}
-1 & 1 & 0 & 0& 0\\
 & & 0 & &  
\end{pmatrix} 
& 
1
&
D_5
&
\hbox{yes (non-singular)}
&
1
\\
\\
3
&
\begin{pmatrix}
0 & -1 & 1 & 0& 0\\
 & & 0 & &  
\end{pmatrix} 
& 
31
&
SL_2 \times SL_5
&
\hbox{yes (non-singular)} 
&
2
\\
\\
4
&
\begin{pmatrix}
0 & 0 & -1 & 1& 0\\
 & & 1 & &  
\end{pmatrix} 
& 
431
&
SL_3 \times SL_2 \times SL_3
&
\hbox{yes (non-singular)}
&
3
\\
5
&
 \begin{pmatrix}
0 & 0 & 0 & 1& 0\\
 & & -1 & &  
\end{pmatrix} 
& 
2431
&
SL_6
&
\hbox{yes (non-singular)}
&
4
\\
\\
6
&
\begin{pmatrix}
0 & 0 & 0 & -1 & 1\\
 & & 1 & &  
\end{pmatrix} 
& 
5431
&
SL_5 \times SL_2
&
\hbox{yes (non-singular)}
&
4
\\
\\
7
&
\begin{pmatrix}
0 & 0 & 1 & -1 & 1\\
 & & -1 & &  
\end{pmatrix} 
& 
52431
&
SL_4 \times SL_2&
\hbox{yes}
&
5
\\
8
&
\begin{pmatrix}
0 & 0 & 0 & 0& -1\\
 & & 1 & &  
\end{pmatrix} 
& 
65431
&
D_5
&
\hbox{yes (non-singular)}
&
5
\\
\\
9
&
\begin{pmatrix}
0 & 1 & -1 & 0 & 1\\
 & & 0 & &  
\end{pmatrix} 
& 
452431
&
SL_3 \times SL_2 \times SL_3
&
\hbox{yes}
&
6
\\
10
&
\begin{pmatrix}
0 & 0 & 1 & 0& -1\\
 & & -1 & &  
\end{pmatrix} 
& 
652431
&
SL_5
&
\hbox{no}
&
6
\\
\\
11
&
\begin{pmatrix}
1 & -1 & 0 & 0& 1\\
 & & 0 & &  
\end{pmatrix} 
& 
3452431
&
SL_2 \times SL_5
&
\hbox{yes}
&
7
\\
12
&
\begin{pmatrix}
0 & 1 & -1 & 1 & -1\\
 & & 0 & &  
\end{pmatrix} 
& 
6452431
&
SL_3 \times SL_2 \times SL_2
&
\hbox{yes}
&
7
\\
\\
13
&
\begin{pmatrix}
1 & -1 & 0 & 1 & -1 \\
 & & 0 & &  
\end{pmatrix} 
& 
63452431
&
SL_2 \times SL_4
&
\hbox{no}
&
8
\\
\\
14
&\begin{pmatrix}
0 & 1 & 0 & -1 & 0\\
 & & 0 & &  
\end{pmatrix} 
& 
54652431
&
SL_5 \times SL_2 
&
\hbox{yes}
&
8
\\
\\
15
&
\begin{pmatrix}
1 & -1 & 1 & -1 & 0\\
 & & 0 & &  
\end{pmatrix} 
& 
563452431
&
SL_2 \times SL_3 \times SL_2
&
\hbox{yes}
&
9
\\
\\
16
&
\begin{pmatrix}
1 & 0 & -1 & 0 & 0\\
 & & 1 & &  
\end{pmatrix} 
& 
4563452431
&
SL_3 \times SL_2 \times SL_3
&
\hbox{yes}
&
10
\\
\\
17
&
\begin{pmatrix}
1 & 0 & 0 & 0& 0\\
 & & -1 & &  
\end{pmatrix} 
& 
24563452431
&
SL_6
&
\hbox{yes}
&
11
\\
\\
18
&
\begin{pmatrix}
-1 & 0 & 0 & 0& 1\\
 & & 0 & &  
\end{pmatrix} 
& 
13452431
&
D_5
&
\hbox{yes}
&
8
\\
\\
19
&
\begin{pmatrix}
-1 & 0 & 0 & 1 & -1\\
 & & 0 & &  
\end{pmatrix} 
& 
163452431
&
D_4
&
\hbox{no}
&
9
\\
\\
20
&
\begin{pmatrix}
-1 & 0 & 1 & -1& 0\\
 & & 0 & &  
\end{pmatrix} 
& 
5613452431
&
SL_3 \times SL_2 \times SL_2
&
\hbox{no}
&
10
\\
\\
21
&
\begin{pmatrix}
-1 & 1 & -1 & 0& 0\\
 & & 1 & &  
\end{pmatrix} 
& 
45613452431
&
SL_2 \times SL_2 \times SL_3
&
\hbox{yes}
&
11
\\
\\
22
&
\begin{pmatrix}
-1 & 1 & 0 & 0& 0\\
 & & -1 & &  
\end{pmatrix} 
& 
245613452431
&
SL_5
&
\hbox{no}
&
12
\\
\\
23
&
\begin{pmatrix}
0 & -1 & 0 & 0& 0\\
 & & 1 & &  
\end{pmatrix} 
& 
345613452431
&
SL_2 \times SL_5
&
\hbox{yes}
&
12
\\
24
&
\begin{pmatrix}
0 & -1 & 1 & 0& 0\\
 & & -1 & &  
\end{pmatrix} 
& 
2345613452431
&
SL_2 \times SL_4
&
\hbox{yes}
&
13
\\
25
&
\begin{pmatrix}
0 & 0 & -1 & 1 & 0\\
 & & 0 & &  
\end{pmatrix} 
& 
42345613452431
&
SL_3 \times SL_2 \times SL_3
&
\hbox{yes}
&
14
\\
26
&
\begin{pmatrix}
0 & 0 & 0 & -1& 1\\
 & & 0 & &  
\end{pmatrix} 
& 
542345613452431
&
SL_5 \times SL_2
&
\hbox{yes}
&
15
\\
27
&
\begin{pmatrix}
0 & 0 & 0 & 0& -1\\
 & & 0 & &  
\end{pmatrix} 
& 
6542345613452431
&
D_5
&
\hbox{yes}
&
16
\end{array}}
$$

\subsection{Weights and symmetries for Schubert varieties in $E_{7}/P_{7}$}
$$
{\renewcommand{\arraystretch}{0.8}
\begin{array}{lllllll}
\mbox{No}
&
\mbox{Weights of $V(\varpi_1)$}
&
\hbox{elt. in }
W/W_{P}
&
\mbox{Group of symmetries}

&
\mbox{Gorenstein}
&
\dim 
\\ 
\hline
\\
1
&
\begin{pmatrix}
0 & 0 & 0 & 0& 0 & 1\\
 & & 0 & &  
\end{pmatrix} 
& 
e
&
E_7
&
\hbox{yes}
&
0
\\
2
&
\begin{pmatrix}
0 & 0 & 0 & 0& 1 & -1\\
 & & 0 & &   
\end{pmatrix} 
& 
7
&
D_6
&
\hbox{yes}
&
1
\\
3
&
\begin{pmatrix}
0 & 0 & 0 & 1& -1 & 0\\
 & & 0 & &    
\end{pmatrix} 
& 
67
&
SL_2 \times D_5
&
\hbox{yes}
&
2
\\
4
&
\begin{pmatrix}
0 & 0 & 1 & -1& 0 & 0\\
 & & 0 & &      
\end{pmatrix} 
& 
567
&
 D_4 \times SL_3
&
\hbox{yes}
&
3
\\
5
&
\begin{pmatrix}
0 & 1 & -1 & 0 & 0 & 0\\
  &   & 1  &   &    
\end{pmatrix} 
& 
4567
&
SL_3 \times SL_2 \times SL_4
&
\hbox{yes}
&
4
\\
6
&
 \begin{pmatrix}
0 & 1 & 0 & 0 & 0 & 0\\
  &   & -1  &   &   
\end{pmatrix} 
& 
24567
&
SL_7
&
\hbox{yes}
&
5
\\
\\
7
&
\begin{pmatrix}
1 & -1 & 0 & 0 & 0 & 0\\
  &    & 1  &   &   
\end{pmatrix} 
& 
34567
&
SL_2 \times SL_6
&
\hbox{yes}
&
5
\\
8
&
\begin{pmatrix}
 -1& 0 & 0  & 0 & 0 & 0 \\
  &   & 1  &   &  
\end{pmatrix} 
& 
134567
&
D_6
&
\hbox{yes}
&
6
\\
\\
9
&
\begin{pmatrix}
 1& -1 & 1 & 0 & 0 & 0\\
  &    & -1  &   &  
\end{pmatrix} 
& 
234567
&
SL_2 \times SL_5
&
\hbox{yes}
&
6
\\
\\
10
&
\begin{pmatrix}
-1 & 0 & 1  & 0 & 0 & 0\\
  &   & -1  &   &  
\end{pmatrix} 
& 
2134567
&
SL_6
&
\hbox{no}
&
7
\\
\\
11
&
\begin{pmatrix}
1 & 1 & -1  & 1 & 0 & 0\\
  &   & 0  &   &  
\end{pmatrix} 
& 
4324567
&
SL_3 \times SL_2 \times SL_4
&
\hbox{yes}
&
7
\\
\\
12
&
\begin{pmatrix}
-1 & 1 & -1  & 1 & 0 & 0\\
  &   & 0  &   &   
\end{pmatrix} 
& 
42134567
&
SL_2 \times SL_2 \times SL_4
&
\hbox{yes}
&
8
\\
\\
13
&
\begin{pmatrix}
 1& 1 & 0  & -1 & 1 & 0\\
  &   & 0  &   &  
\end{pmatrix} 
& 
5434567
&
D_4 \times SL_3
&
\hbox{yes}
&
8
\\
14
&
\begin{pmatrix}
0 & -1 & 0  & 1 & 0 & 0\\
  &   & 0  &   &   
\end{pmatrix} 
& 
342134567
&
SL_2 \times SL_6
&
\hbox{no}
&
9
\\
15
&
\begin{pmatrix}
-1 & 1 & 0  & -1 & 1 & 0\\
  &   & 0  &   &   
\end{pmatrix} 
& 
542134567
&
SL_4 \times SL_3
&
\hbox{yes}
&
9
\\
\\
16
&
\begin{pmatrix}
 1& 1 & 0  & 0 & -1 & 1\\
  &   & 0  &   &  
\end{pmatrix} 
& 
654234567
&
D_5 \times SL_2
&
\hbox{no}
&
9
\\
\\
17
&
\begin{pmatrix}
0 & -1 & 1 & -1 & 1 & 0\\
  &   & 0  &   &   
\end{pmatrix} 
& 
6342134567
&
SL_2 \times SL_3 \times SL_3
&
\hbox{yes}
&
10
\\
\\
18
&
\begin{pmatrix}
-1 & 1 & 0  & 0 & -1 & 1\\
  &   & 0  &   &   
\end{pmatrix} 
& 
6542134567
&
D_4 \times SL_2
&
\hbox{no}
&
10
\\
\\
19
&
\begin{pmatrix}
0 & 0 & -1 & 0 & 1 & 0\\
  &   & 1  &   &   
\end{pmatrix} 
& 
43542134567
&
SL_3 \times SL_2 \times SL_4
&
\hbox{no}
&
11
\\
20
&
\begin{pmatrix}
0 & -1 & 1 & 0 & -1 & 1\\
  &   & 0  &   &   
\end{pmatrix} 
& 
65342134567
&
SL_2 \times SL_4 \times SL_2
&
\hbox{no}
&
11
\\
\\
21
&
\begin{pmatrix}
0 & 0 & 0 & 0 & 1 & 0\\
  &   & -1  &   &   
\end{pmatrix} 
& 
245342134567
&
SL_7
&
\hbox{no}
&
12
\\
\\
22
&
\begin{pmatrix}
0 & 0 & -1 & 1 & -1 & 1\\
  &   & 1  &   &   
\end{pmatrix} 
& 
645342134567
&
SL_3 \times SL_2 \times SL_2 \times SL_2
&
\hbox{no}
&
12
\\
\\
23
&
\begin{pmatrix}
0 & 0 & 0 & 1 & -1 & 1\\
  &   & -1  &   &   
\end{pmatrix} 
& 
2645342134567
&
SL_5 \times SL_2
&
\hbox{no}
&
13
\\
\\
24
&
\begin{pmatrix}
0 & 0 & 0 & -1 & 0 & 1\\
  &   & 1  &   &   
\end{pmatrix} 
& 
5645342134567
&
D_4 \times SL_3
&
\hbox{yes}
&
13
\\
25
&
\begin{pmatrix}
0 & 0 & 1 & -1 & 0 & 1\\
  &   & -1  &   &   
\end{pmatrix} 
& 
25645342134567
&
SL_4 \times SL_3
&
\hbox{yes}
&
14
\\
26
&
\begin{pmatrix}
0 & 1 & -1 & 0 & 0 & 1\\
  &   & 0  &   &   
\end{pmatrix} 
& 
425645342134567
&
SL_3 \times SL_2 \times SL_4
&
\hbox{yes}
&
15
\\
27
&
\begin{pmatrix}
1 & -1 & 0 & 0 & 0 & 1\\
  &   & 0  &   &   
\end{pmatrix} 
& 
3425645342134567
&
SL_2 \times SL_6
&
\hbox{yes}
&
16
\\
28
&
\begin{pmatrix}
-1 & 0 & 0 & 0 & 0 & 1\\
  &   & 0  &   &   
\end{pmatrix} 
& 
13425645342134567
&
D_6
&
\hbox{yes}
&
17
\end{array}}
$$

$$
{\renewcommand{\arraystretch}{0.8}
\begin{array}{lllllll}
\mbox{No}
&
\mbox{Weights of $V(\varpi_1)$}
&
\hbox{elt. in }
W/W_{P}
&
\mbox{Group of symmetries}

&
\mbox{Gorenstein}
&
\dim 
\\ 
\hline
\\
29
&
\begin{pmatrix}
 1& 1 & 0  & 0 & 0 & -1\\
  &   & 0  &   &  
\end{pmatrix} 
& 
7654234567
&
E_6
&
\hbox{yes (non-singular)}
&
10
\\
30
&
\begin{pmatrix}
 -1& 2 & 0  & 0 & 0 & -1\\
  &   & 0  &   &  
\end{pmatrix} 
& 
17654234567
&
D_5
&
\hbox{yes (non-singular)}
&
11
\\
31
&
\begin{pmatrix}
 1& -2 & 2  & 0 & 0 & -1\\
  &   & 0  &   &  
\end{pmatrix} 
& 
317654234567
&
SL_2 \times SL_5
&
\hbox{no}
&
12
\\
32
&
\begin{pmatrix}
 1& 0 & -2  & 2 & 0 & -1\\
  &   & 2  &   &  
\end{pmatrix} 
& 
4317654234567
&
SL_3 \times SL_2 \times SL_3
&
\hbox{no}
&
13
\\
33
&
\begin{pmatrix}
 1& 0 & 0  & 2 & 0 & -1\\
  &   & -2  &   &  
\end{pmatrix} 
& 
24317654234567
&
SL_6 
&
\hbox{yes} 
&
14
\\
34
&
\begin{pmatrix}
 1& 0 & 0  & -2 & 2 & -1\\
  &   & 2  &   &  
\end{pmatrix} 
& 
54317654234567
&
D_5 \times SL_2 
&
\hbox{yes} 
&
14
\\
35
&
\begin{pmatrix}
 1& 0 & 2  & -2 & 2 & -1\\
  &   & -2  &   &  
\end{pmatrix} 
& 
524317654234567
&
SL_4 \times SL_2 
&
\hbox{yes} 
&
15
\\
36
&
\begin{pmatrix}
 1& 0 & 0  & 0 & -2 & 1\\
  &   & 2  &   &  
\end{pmatrix} 
& 
654317654234567
&
D_5 \times SL_2 
&
\hbox{yes} 
&
15
\\
37
&
\begin{pmatrix}
 1& 2 & -2  & 0 & 2 & -1\\
  &   & 0  &   &  
\end{pmatrix} 
& 
4524317654234567
&
SL_3 \times SL_2 \times SL_3
&
\hbox{no} 
&
16
\\
38
&
\begin{pmatrix}
 1& 0 & 2  & 0 & -2 & 1\\
  &   & -2  &   &  
\end{pmatrix} 
& 
2654317654234567
&
SL_5 \times SL_2 
&
\hbox{yes} 
&
16
\\
39
&
\begin{pmatrix}
3& -2 & 0  & 0 & 2 & -1\\
  &   & 0  &   &  
\end{pmatrix} 
& 
34524317654234567
&
SL_2 \times SL_5 
&
\hbox{no} 
&
17
\\
40
&
\begin{pmatrix}
 1& 2 & -2  & 2 & -2 & 1\\
  &   & 0  &   &  
\end{pmatrix} 
& 
42654317654234567
&
SL_3 \times SL_2 \times SL_2 \times SL_2
&
\hbox{yes} 
&
17
\\
41
&
\begin{pmatrix}
-3& 1 & 0  & 0 & 2 & -1\\
  &   & 0  &   &  
\end{pmatrix} 
& 
134524317654234567
&
SL_3 \times D_5
&
\hbox{no} 
&
18
\\
42
&
\begin{pmatrix}
 3& -2 & 0  & 2 & -2 & 1\\
  &   & 0  &   &  
\end{pmatrix} 
& 
342654317654234567
&
SL_2 \times SL_4 \times SL_2
&
\hbox{yes} 
&
18
\\
43
&
\begin{pmatrix}
 1& 2 & 0  & -2 & 0 & 1\\
  &   & 0  &   &  
\end{pmatrix} 
& 
542654317654234567
&
SL_5 \times SL_3
&
\hbox{yes} 
&
18
\\
44
&
\begin{pmatrix}
-3& 1 & 0  & 2 & -2 & 1\\
  &   & 0  &   &  
\end{pmatrix} 
& 
6134524317654234567
&
D_4 \times SL_2
&
\hbox{no} 
&
19
\\
45
&
\begin{pmatrix}
 3& -2 & 2  & -2 & 0 & 1\\
  &   & 0  &   &  
\end{pmatrix} 
& 
3542654317654234567
&
SL_2 \times SL_3 \times SL_3
&
\hbox{yes} 
&
19
\\
46
&
\begin{pmatrix}
-3& 1 & 2  & -2 & 0 & 1\\
  &   & 0  &   &  
\end{pmatrix} 
& 
56134524317654234567
&
SL_4 \times SL_3
&
\hbox{no} 
&
20
\\
47
&
\begin{pmatrix}
 3& 0 & -2  & 0 & 0 & 1\\
  &   & 2  &   &  
\end{pmatrix} 
& 
43542654317654234567
&
SL_3 \times SL_2 \times SL_4
&
\hbox{yes} 
&
20
\\
48
&
\begin{pmatrix}
-1 & 1 & -1 & 0 & 0 & 0 \\
   &  & 1 & &  
\end{pmatrix} 
& 
145673425645341234567
&
SL_2 \times SL_2 \times SL_4
&
\hbox{yes}
&
21
\\
\\
49
&
\begin{pmatrix}
1 & 0 & 0 & 0 & 0 & 0 \\
   &  & -1 & &  
\end{pmatrix} 
& 
245673425645341234567
&
SL_7
&
\hbox{yes}
&
21
\\
50
&
\begin{pmatrix}
-1 & 1 & 0 & 0 & 0 & 0 \\
   &  & -1 & &  
\end{pmatrix} 
& 
1245673425645341234567
&
SL_6
&
\hbox{no}
&
22
\\
\\
51
&
\begin{pmatrix}
0 & -1 & 0 & 0 & 0 & 0\\
   &  & 1 & &  
\end{pmatrix} 
& 
3145673425645341234567
&
SL_2 \times SL_6
&
\hbox{yes}
&
22
\\
\\
52
&
\begin{pmatrix}
0 & -1 & 1 & 0 & 0 & 0 \\
   &  & -1 & &  
\end{pmatrix} 
& 
31245673425645341234567
&
SL_2 \times SL_5
&
\hbox{yes}
&
23
\\
53
&
\begin{pmatrix}
0 & 0 & -1 & 1 & 0 & 0\\
   &  & 0 & &  
\end{pmatrix} 
& 
431245673425645341234567
&
SL_3 \times SL_2 \times SL_4
&
\hbox{yes}
&
24
\\
54
&
\begin{pmatrix}
0 & 0 & 0 & -1 & 1& 0 \\
   &  & 0 & &  
\end{pmatrix} 
& 
5431245673425645341234567
&
SL_5 \times SL_3
&
\hbox{yes}
&
25
\\
55
&
\begin{pmatrix}
0 & 0 & 0 & 0 & -1& 1\\
   &  & 0 & &  
\end{pmatrix} 
& 
65431245673425645341234567
&
D_5 \times SL_2
&
\hbox{yes}
&
26
\\
56
&
\begin{pmatrix}
0 & 0 & 0 & 0 & 0 & -1\\
  &   & 0 &   &  
\end{pmatrix} 
&
765431245673425645341234567 
&
E_6
&
\hbox{yes}
&
27
\end{array}}
$$

\section{Cubic and derivatives}
\subsection{Cubic in different labellings}
\label{app:cubic}
\begin{align}
\label{cubice6d5}
Q_{E_6|D_5} = & \quad 
x_{1}x_{18}x_{27} - x_{1}x_{19}x_{26} + x_{1}x_{20}x_{25} -  x_{1}x_{21}x_{24}
+  x_{1}x_{22}x_{23} -x_{2}x_{11}x_{27} + x_{2}x_{13}x_{26 }- x_{2}x_{15}x_{25}  
+ x_{2}x_{16}x_{24} 
\nonumber \\ 
& - x_{2}x_{17}x_{23} +x_{3}x_{9}x_{27} - x_{3}x_{12}x_{26} +x_{3}x_{14}x_{25} - x_{3}x_{16}x_{22} + x_{3}x_{17}x_{21} -x_{4}x_{7}x_{27} +x_{4}x_{10}x_{26} -x_{4}x_{14}x_{24} 
\nonumber \\
&
+x_{4}x_{15}x_{22} -x_{4}x_{17}x_{20} +x_{5}x_{6}x_{27} -x_{5}x_{8}x_{26} + x_{5}x_{14}x_{23} -x_{5}x_{15}x_{21} +x_{5}x_{16}x_{20} -x_{6}x_{10}x_{25} + x_{6}x_{12}x_{24} 
\nonumber \\
&
- x_{6}x_{13}x_{22} + x_{6}x_{17}x_{19} + x_{7}x_{8}x_{25} -x_{7}x_{12}x_{23} + x_{7}x_{13}x_{21} - x_{7}x_{16}x_{19} - x_{8}x_{9}x_{24} + x_{8}x_{11}x_{22} -x_{8}x_{17}x_{18} 
\nonumber \\
&
+ x_{9}x_{10}x_{23}  - x_{9}x_{13}x_{20} + x_{9}x_{15}x_{19} - x_{10}x_{11}x_{21} +x_{10}x_{16}x_{18} + x_{11}x_{12}x_{20} - x_{11}x_{14}x_{19} - x_{12}x_{15}x_{18} + x_{13}x_{14}x_{18} \nonumber
\end{align}

$
\begin{array}{ll}
\label{derive6d5}
f_{1} = z_{1}\overline{z}_{1} - z_{2}\overline{z}_{2} + z_{3}\overline{z}_{3} - z_{4}\overline{z}_{4} + z_{5}\overline{z}_{5} 
&
f_{2} = -y_{2345}\overline{z}_{1} + y_{1345}\overline{z}_{2} - y_{1245}\overline{z}_{3} + y_{1235}\overline{z}_{4} - y_{1234}\overline{z}_{5}
\\
f_{3} = y_{23}\overline{z}_{1} - y_{13}\overline{z}_{2} + y_{12}\overline{z}_{3} - y_{1235}z_{5} + y_{1234}z_{4}
&
f_{4} = -y_{24}\overline{z}_{1} + y_{25}\overline{z}_{2} - y_{12}\overline{z}_{4} + y_{1245}z_{5} - y_{1234}z_{3}
\\
f_{5} = y_{25}\overline{z}_{1} - y_{15}\overline{z}_{2} + y_{12}\overline{z}_{5} - y_{1245}z_{4} + y_{1235}z_{3}
&
f_{6} = y_{34}\overline{z}_{1} - y_{25}\overline{z}_{3} + y_{13}\overline{z}_{4} - y_{1345}z_{5} + y_{1234}z_{2}
\\
f_{7} = -y_{13}\overline{z}_{1} + y_{15}\overline{z}_{3} - y_{13}\overline{z}_{5} + y_{1345}z_{4} - y_{1235}z_{2}
&
f_{8} = -y_{34}\overline{z}_{2} + y_{24}\overline{z}_{3} - y_{23}\overline{z}_{4} + y_{2345}z_{5} - y_{1234}z_{1}
\\
f_{9} = y_{45}\overline{z}_{1} - y_{15}\overline{z}_{4} + y_{25}\overline{z}_{5} - y_{1345}z_{3} + y_{1245}z_{2}
&
f_{10} = y_{13}\overline{z}_{2} - y_{25}\overline{z}_{3} + y_{23}\overline{z}_{5} - y_{2345}z_{4} + y_{1235}z_{1}
\\
f_{11} = -y_{0}\overline{z}_{1} + y_{15}z_{5} - y_{25}z_{4} + y_{13}z_{3} - y_{12}z_{2}
&
f_{12} = -y_{45}\overline{z}_{2} + y_{25}\overline{z}_{4} - y_{24}\overline{z}_{5} + y_{2345}z_{3} - y_{1245}z_{1}
\\
f_{13} = y_{0}\overline{z}_{2} - y_{25}z_{5} + y_{24}z_{4} - y_{23}z_{3} + y_{12}z_{1}
&
f_{14} = y_{45}\overline{z}_{3} - y_{13}\overline{z}_{4} + y_{34}\overline{z}_{5} - y_{2345}z_{2} + y_{1345}z_{1}
\\
f_{15} = -y_{0}\overline{z}_{3} + y_{13}z_{5} - y_{34}z_{4} + y_{23}z_{2} - y_{13}z_{1}
&
f_{16} = y_{0}\overline{z}_{4} - y_{45}z_{5} + y_{34}z_{3} - y_{24}z_{2} + y_{25}z_{1}
\\
f_{17} = -y_{0}\overline{z}_{5} + y_{45}z_{4} - y_{13}z_{3} + y_{25}z_{2} - y_{15}z_{1}
&
f_{18} = x\overline{z}_{1} - y_{15}y_{1234} + y_{25}y_{1235} - y_{13}y_{1245} + y_{1345}y_{12}
\\
f_{19} = -x\overline{z}_{2} + y_{25}y_{1234} - y_{24}y_{1235} + y_{23}y_{1245} - y_{2345}y_{12}
&
f_{20} = x\overline{z}_{3} - y_{13}y_{1234} + y_{34}y_{1235} - y_{23}y_{1345} + y_{2345}y_{13}
\\
f_{21} = -x\overline{z}_{4} + y_{45}y_{1234} - y_{34}y_{1245} + y_{24}y_{1345} - y_{25}y_{2345}
&
f_{22} = x\overline{z}_{5} - y_{45}y_{1235} + y_{13}y_{1245} - y_{25}y_{1345} + y_{15}y_{2345}
\\
f_{23} = xz_{5} - y_{0}y_{1234} + y_{34}y_{12} - y_{24}y_{13} + y_{23}y_{25}
&
f_{24} = -xz_{4} + y_{0}y_{1235} - y_{13}y_{12} + y_{25}y_{13} - y_{15}y_{23}
\\
f_{25} = xz_{3} - y_{0}y_{1245} + y_{45}y_{12} - y_{25}y_{25} + y_{24}y_{15}
&
f_{26} = -xz_{2} + y_{0}y_{1345} - y_{45}y_{13} + y_{13}y_{25} - y_{34}y_{15}
\\
f_{27} = xz_{1} - y_{0}y_{2345} + y_{45}y_{23} - y_{13}y_{24} + y_{34}y_{25}
\end{array}
$

\subsection{Derivatives of the cubic in different labellings}
Cubic in $E_7 \downarrow E_6$-numbering:\\

\begin{align*}
Q_{E7E6} = & 
x_{13}x_{14}x_{15} - x_{11}x_{15}x_{16} - x_{12}x_{13}x_{17} + 
x_{10}x_{16}x_{17} + x_{9}x_{15}x_{18} - x_{8}x_{17}x_{18}  + 
x_{11}x_{12}x_{19} - x_{10}x_{14}x_{19} + x_{5}x_{18}x_{19}
\\ 
&
- x_{7}x_{15}x_{20} + x_{6}x_{17}x_{20} - x_{4}x_{19}x_{20} 
- x_{9}x_{12}x_{21} + x_{8}x_{14}x_{21} - x_{5}x_{16}x_{21} 
+ x_{3}x_{20}x_{21} + x_{7}x_{12}x_{22} - x_{6}x_{14}x_{22} 
\\ 
&
+ x_{4}x_{16}x_{22} - x_{3}x_{18}x_{22} + x_{9}x_{10}x_{23} 
- x_{8}x_{11}x_{23} + x_{5}x_{13}x_{23} - x_{2}x_{20}x_{23} 
+ x_{1}x_{22}x_{23} - x_{7}x_{10}x_{24} + x_{6}x_{11}x_{24} 
\\ 
&
- x_{4}x_{13}x_{24} + x_{2}x_{18}x_{24} - x_{1}x_{21}x_{24} 
+ x_{7}x_{8}x_{25} - x_{6}x_{9}x_{25} + x_{3}x_{13}x_{25} 
- x_{2}x_{16}x_{25} + x_{1}x_{19}x_{25} - x_{5}x_{7}x_{26} 
\\ 
&
+ x_{4}x_{9}x_{26} - x_{3}x_{11}x_{26} + x_{2}x_{14}x_{26} 
- x_{1}x_{17}x_{26} + x_{5}x_{6}x_{27} - x_{4}x_{8} x_{27} 
+ x_{3}x_{10}x_{27} - x_{2}x_{12}x_{27} + x_{1}x_{15}x_{27}
\end{align*}

Denote by $f_i$ the derivative of $Q_{E7E6}$ with respect to $x_i$ for $i=1, \ldots, 27$:

$
\begin{array}{ll}
\label{derive7e6}
 f_{1} = {x}_{22} {x}_{23}-{x}_{21} {x}_{24}+{x}_{19} {x}_{25}-{x}_{17} {x}_{26}+{x}_{15} {x}_{27}  & 
 f_{2} = -{x}_{20} {x}_{23}+{x}_{18} {x}_{24}-{x}_{16} {x}_{25}+{x}_{14} {x}_{26}-{x}_{12} {x}_{27}  \\ 
 f_{3} = {x}_{20} {x}_{21}-{x}_{18} {x}_{22}+{x}_{13} {x}_{25}-{x}_{11} {x}_{26}+{x}_{10} {x}_{27}  & 
 f_{4} = -{x}_{19} {x}_{20}+{x}_{16} {x}_{22}-{x}_{13} {x}_{24}+{x}_{9} {x}_{26}-{x}_{8} {x}_{27}  \\ 
 f_{5} = {x}_{18} {x}_{19}-{x}_{16} {x}_{21}+{x}_{13} {x}_{23}-{x}_{7} {x}_{26}+{x}_{6} {x}_{27}  &
 f_{6} = {x}_{17} {x}_{20}-{x}_{14} {x}_{22}+{x}_{11} {x}_{24}-{x}_{9} {x}_{25}+{x}_{5} {x}_{27}  \\ 
 f_{7} = -{x}_{15} {x}_{20}+{x}_{12} {x}_{22}-{x}_{10} {x}_{24}+{x}_{8} {x}_{25}-{x}_{5} {x}_{26}  & 
 f_{8} = -{x}_{17} {x}_{18}+{x}_{14} {x}_{21}-{x}_{11} {x}_{23}+{x}_{7} {x}_{25}-{x}_{4} {x}_{27}  \\ 
 f_{9} = {x}_{15} {x}_{18}-{x}_{12} {x}_{21}+{x}_{10} {x}_{23}-{x}_{6} {x}_{25}+{x}_{4} {x}_{26}  & 
 f_{10} = {x}_{16} {x}_{17}-{x}_{14} {x}_{19}+{x}_{9} {x}_{23}-{x}_{7} {x}_{24}+{x}_{3} {x}_{27}  \\ 
 f_{11} = -{x}_{15} {x}_{16}+{x}_{12} {x}_{19}-{x}_{8} {x}_{23}+{x}_{6} {x}_{24}-{x}_{3} {x}_{26}  & 
 f_{12} = -{x}_{13} {x}_{17}+{x}_{11} {x}_{19}-{x}_{9} {x}_{21}+{x}_{7} {x}_{22}-{x}_{2} {x}_{27}  \\ 
 f_{13} = {x}_{14} {x}_{15}-{x}_{12} {x}_{17}+{x}_{5} {x}_{23}-{x}_{4} {x}_{24}+{x}_{3} {x}_{25}  & 
 f_{14} = {x}_{13} {x}_{15}-{x}_{10} {x}_{19}+{x}_{8} {x}_{21}-{x}_{6} {x}_{22}+{x}_{2} {x}_{26}  \\ 
 f_{15} = {x}_{13} {x}_{14}-{x}_{11} {x}_{16}+{x}_{9} {x}_{18}-{x}_{7} {x}_{20}+{x}_{1} {x}_{27}  & 
 f_{16} = -{x}_{11} {x}_{15}+{x}_{10} {x}_{17}-{x}_{5} {x}_{21}+{x}_{4} {x}_{22}-{x}_{2} {x}_{25}  \\ 
 f_{17} = -{x}_{12} {x}_{13}+{x}_{10} {x}_{16}-{x}_{8} {x}_{18}+{x}_{6} {x}_{20}-{x}_{1} {x}_{26}  & 
 f_{18} = {x}_{9} {x}_{15}-{x}_{8} {x}_{17}+{x}_{5} {x}_{19}-{x}_{3} {x}_{22}+{x}_{2} {x}_{24}  \\ 
 f_{19} = {x}_{11} {x}_{12}-{x}_{10} {x}_{14}+{x}_{5} {x}_{18}-{x}_{4} {x}_{20}+{x}_{1} {x}_{25}  & 
 f_{20} = -{x}_{7} {x}_{15}+{x}_{6} {x}_{17}-{x}_{4} {x}_{19}+{x}_{3} {x}_{21}-{x}_{2} {x}_{23}  \\ 
 f_{21} = -{x}_{9} {x}_{12}+{x}_{8} {x}_{14}-{x}_{5} {x}_{16}+{x}_{3} {x}_{20}-{x}_{1} {x}_{24}  & 
 f_{22} = {x}_{7} {x}_{12}-{x}_{6} {x}_{14}+{x}_{4} {x}_{16}-{x}_{3} {x}_{18}+{x}_{1} {x}_{23}  \\ 
 f_{23} = {x}_{9} {x}_{10}-{x}_{8} {x}_{11}+{x}_{5} {x}_{13}-{x}_{2} {x}_{20}+{x}_{1} {x}_{22}  & 
 f_{24} = -{x}_{7} {x}_{10}+{x}_{6} {x}_{11}-{x}_{4} {x}_{13}+{x}_{2} {x}_{18}-{x}_{1} {x}_{21}  \\ 
 f_{25} = {x}_{7} {x}_{8}-{x}_{6} {x}_{9}+{x}_{3} {x}_{13}-{x}_{2} {x}_{16}+{x}_{1} {x}_{19}  &
 f_{26} = -{x}_{5} {x}_{7}+{x}_{4} {x}_{9}-{x}_{3} {x}_{11}+{x}_{2} {x}_{14}-{x}_{1} {x}_{17}  \\ 
 f_{27} = {x}_{5} {x}_{6}-{x}_{4} {x}_{8}+{x}_{3} {x}_{10}-{x}_{2} {x}_{12}+{x}_{1} {x}_{15}  & \\
\end{array}
$

\section*{Acknowledgements}

J.W. was supported by  the grants MAESTRO NCN - UMO-2019/34/A/ST1/00263 - Research in Commutative Algebra and Representation Theory and NAWA POWROTY - PPN/PPO/2018/1/00013/U/00001 - Applications of Lie algebras to Commutative Algebra. S.A.F. was supported by the grant NAWA POWROTY - PPN/PPO/2018/1/00013/U/00001 - Applications of Lie algebras to Commutative Algebra and J.T. was supported by the grant MAESTRO NCN - UMO-2019/34/A/ST1/00263 - Research in Commutative Algebra and Representation Theory. J.T. was also partially supported by  {\it Narodowe Centrum Nauki}, grant number 2017/26/A/ST1/00189. A lot of progress toward the completion of this manuscript and the ideas behind it was made during the ICERM Meeting \textit{Free Resolutions and Representation Theory}, held online on August 3rd-7th 2020. We would like to thank Lars Christensen, Ela Celikbas, David Eisenbud, Lorenzo Guerreri,  Witold Kra\'skiewicz, Craig Huneke, Andrew Kustin, Jai Laxmi,  Kyu-Hwan Lee, Pedro Macias Marques, Claudia Polini, Bernd Ulrich and  Oana Veliche for fruitful discussions.

\newpage
\restoregeometry
\bibliographystyle{abbrv}

\end{document}